\documentclass[a4paper,11pt]{amsart}
\usepackage{amssymb}
\newtheorem{lemma}{Lemma}
\newtheorem{prop}{Proposition}
\newtheorem{thm}{Theorem}
\newtheorem*{result}{Theorem}

\theoremstyle{definition}
\newtheorem{defn}{Definition}
\theoremstyle{remark}
\newtheorem{rem}{Remark}
\newtheorem{ex}{Example}
\newcounter{numl}

\newcommand{\labelnuml}{\textup{(\roman{numl})}}
\newenvironment{numlist}{\begin{list}{\labelnuml}%
{\usecounter{numl}\setlength{\leftmargin}{0pt}%
\setlength{\itemindent}{2\parindent}%
\setlength{\itemsep}{\smallskipamount}\def
\makelabel ##1{\hss \llap {\upshape ##1}}}}{\end{list}}

\newenvironment{numlproof}{\begin{proof}\begin{numlist}}
{\qed \end{numlist}\begingroup\let\qed\relax\end{proof}\endgroup}
\newcommand{\R}{{\mathbb R}}
\newcommand{\C}{{\mathbb C}}

\newcommand{\cF}{{\mathcal F}}
\newcommand{\cL}{{\mathcal L}}
\newcommand{\Ric}{\mathit{Ric}}
\newcommand{\Scal}{\mathit{Scal}}
\newcommand{\Id}{\mathit{Id}}
\newcommand{\sym}{\mathop{\mathrm{sym}}\nolimits}
\newcommand{\trace}{\mathop{\mathrm{tr}}\nolimits}
\newcommand{\vol}{\mathop{\mathrm{vol}}\nolimits}
\newcommand{\pfaff}{\mathop{\mathrm{pf}}\nolimits}
\newcommand{\grad}{\mathop{\mathrm{grad}}\nolimits}
\newcommand{\g}[1]{\langle #1 \rangle}
\newcommand{\G}[1]{\bigl\langle #1 \bigr\rangle}
\newcommand{\symprod}{\mathbin{\raise1pt\hbox{$\scriptstyle\bigcirc$}}}
\newcommand{\eps}{\varepsilon}
\newcommand{\restr}[1]{|_{#1}^{\vphantom x}}
\newcommand{\nrho}{\tilde\rho}
\newcommand{\pF}{F}
\newcommand{\oF}{\Theta}
\newcommand{\Mp}{p}
\newcommand{\Mpc}{p_{\mathrm{c}}}
\newcommand{\Mpn}{p_{\mathrm{nc}}}
\newcommand{\CP}{F_{\mathrm c}}
\newcommand{\Fm}{F_{\mathrm{m}}}
\newcommand{\WK}{W^{\mathcal K}}
\begin{document}
\title[Hamiltonian 2-forms in K{\smash{\"a}}hler geometry, I]
{Hamiltonian 2-forms in K{\smash{\"a}}hler geometry,\\ I General theory}
\author[V. Apostolov]{Vestislav Apostolov}
\author[D. Calderbank]{David M. J. Calderbank}
\author[P. Gauduchon]{Paul Gauduchon}
\thanks{The first author was supported in part by FCAR grant NC-7264,
and by NSERC grant OGP0023879, the second author by the Leverhulme Trust
and the William Gordon Seggie Brown Trust. All three authors are members
of EDGE, Research Training Network HPRN-CT-2000-00101, supported by the
European Human Potential Programme.}
\address{Vestislav Apostolov \\ D{\'e}partement de Math{\'e}matiques\\
UQAM\\ C.P. 8888 \\ Succ. Centre-ville \\ Montr{\'e}al (Qu{\'e}bec) \\
H3C 3P8 \\ Canada}
\email{apostolo@math.uqam.ca}
\address{David M. J. Calderbank \\ School of Mathematics \\
  University of Edinburgh\\ King's Buildings\\ Mayfield Road\\
  Edinburgh EH9 3JZ\\ Scotland}
\email{davidmjc@maths.ed.ac.uk}
\address{Paul Gauduchon \\ Centre de Math\'ematiques\\
{E}cole Polytechnique \\ UMR 7640 du CNRS
\\ 91128 Palaiseau \\ France}
\email{pg@math.polytechnique.fr}
\date{November 2002}
\begin{abstract} We introduce the notion of a hamiltonian $2$-form on a
K\"ahler manifold and obtain a complete local classification. This notion
appears to play a pivotal role in several aspects of K\"ahler geometry. In
particular, on any K\"ahler manifold with co-closed Bochner tensor, the
(suitably normalized) Ricci form is hamiltonian, and this leads to an explicit
description of these K\"ahler metrics, which we call {\it weakly
Bochner-flat}.  Hamiltonian $2$-forms also arise on conformally Einstein
K\"ahler manifolds and provide an Ansatz for extremal K\"ahler metrics
unifying and extending many previous constructions.
\end{abstract}
\maketitle
\vspace{-4mm}

In a previous paper~\cite{ACG0}, while investigating K\"ahler $4$-manifolds
whose antiselfdual Weyl tensor is co-closed, we happened upon a remarkable
linear differential equation for $(1,1)$-forms $\phi$ on a K\"ahler manifold.
This equation states (in any dimension)
\begin{equation} \label{eq:ham1}
\nabla _X \phi = \frac12 (d \trace_\omega\phi
\wedge g(JX,\cdot) - J d \trace_\omega \phi \wedge g(X,\cdot))
\end{equation}
for all vector fields $X$, where $(g,J,\omega)$ is the K\"ahler structure
with Levi-Civita connection $\nabla$. A {\it hamiltonian $2$-form} is a
(nontrivial) solution $\phi$ of~\eqref{eq:ham1}.

Hamiltonian $2$-forms underpin many explicit constructions in K\"ahler
geometry. They arise in particular on Bochner-flat K\"ahler manifolds and on
K\"ahler manifolds (of dimension greater than four) which are conformally
Einstein, both of which have been classified recently, respectively by
Bryant~\cite{bryant}, and Derdzi\'nski and Maschler~\cite{derd-masch}.  In
this paper we obtain an explicit local classification of all K\"ahler metrics
with a hamiltonian $2$-form, which provides a unifying framework for these
works, and at the same time extends Bryant's local classification to the much
larger class of K\"ahler manifolds with co-closed Bochner tensor, called {\it
weakly Bochner-flat}.

The key feature of hamiltonian $2$-forms $\phi$ on K\"ahler $2m$-manifolds
$M$---and the reason for the name---is that if $\sigma_1,\ldots\sigma_m$ are
the elementary symmetric functions of the $m$ eigenvalues of $\phi$ (viewed as
a hermitian operator via the K\"ahler form $\omega$), then the hamiltonian
vector fields $K_r=J\grad_g \sigma_r$ are Killing. Further, the Poisson
brackets $\{\sigma_r, \sigma_s\}$ are all zero, so that the vector fields
$K_1,\ldots K_m$ commute.

If $K_1,\ldots K_m$ are linearly independent, then the K\"ahler metric is
toric. However, not every toric K\"ahler metric arises in this way: the
hamiltonian property also implies that the eigenvalues of $\phi$ have
orthogonal gradients. We say that a toric manifold is {\it orthotoric} if
there is a momentum map $(\sigma_1,\ldots\sigma_m)$ for the torus action (with
respect to some basis of the Lie algebra) such that the gradients of the roots
of the polynomial $\sum_{r=0}^m (-1)^r \sigma_r t^{m-r}$ are orthogonal---here
$\sigma_0=1$.

Of course $K_1,\ldots K_m$ need not be independent; then on an open set where
the span is $\ell$-dimensional, there is a local hamiltonian $\ell$-torus
action by isometries, so the K\"ahler metric on $M$ may be described (locally)
by the Pedersen--Poon construction~\cite{ped-poon}, as a fibration, with
$2\ell$-dimensional toric fibres, over a $2(m-\ell)$-dimensional complex
manifold $S$ equipped with a family of K\"ahler quotient metrics parameterized
by the momentum map of the local $\ell$-torus action.

The hamiltonian property of $\phi$ has further implications for the geometry
of this fibration and of the base $S$.  We say that a hamiltonian $\ell$-torus
action is {\it rigid} if the metric on the orbits depends only on the momentum
map, and that the base $S$ is {\it semisimple} if the K\"ahler quotient
metrics are simultaneously diagonalizable and have common Levi-Civita
connection. The rigidity condition has its origins in work of Calabi on
K\"ahler metrics on holomorphic bundles~\cite{calabi0,calabi1} and has a
number of formulations: it means, for instance, that the local fibration of
$M$ over $S$ is totally geodesic, or equivalently, that $M$ is associated,
locally, to a principal $\ell$-torus bundle with connection over $S$. When
$\ell=1$ semisimplicity is closely related to the `$\sigma$-constancy' of
Hwang--Singer~\cite{hwang-singer}.  Both rigidity and semisimplicity have
explicit descriptions as special cases of the Pedersen--Poon construction.

Our main result shows that any K\"ahler manifold with a hamiltonian $2$-form
induces a (local) semisimple rigid $\ell$-torus action $(\ell\leq m)$ with
orthotoric fibres, and that conversely, such an explicit K\"ahler metric has a
hamiltonian $2$-form. From this we deduce a local classification of weakly
Bochner-flat K\"ahler metrics, rederiving in particular Bryant's
classification of Bochner-flat K\"ahler metrics. We also obtain a new proof of
the Derdzi\'nski--Maschler classification of conformally Einstein K\"ahler
metrics in higher dimensions, and an Ansatz for extremal K\"ahler
metrics---including all constant scalar curvature and K\"ahler--Einstein
metrics with a hamiltonian $2$-form---which unifies and extends many
constructions in the literature.

The structure of the paper is as follows. In section~\ref{s:ckm}, after
reviewing some background material, we explain how equation~\eqref{eq:ham1}
arises on weakly Bochner-flat K\"ahler manifolds and on conformally Einstein
K\"ahler manifolds. Thus motivated, we begin the study of hamiltonian
$2$-forms in section~\ref{s:ham}, where we derive the existence of the
hamiltonian Killing vector fields, and show that the equation for hamiltonian
$2$-forms is an overdetermined equation of finite type which is completely
integrable on manifolds of constant holomorphic sectional curvature.

In section~\ref{s:hta} we study (isometric) hamiltonian torus actions in
general. This section is almost entirely independent of the first two,
although our analysis is motivated by the special properties of hamiltonian
$2$-forms. We first show that the Pedersen--Poon construction~\cite{ped-poon}
has a natural and essentially coordinate-free description in terms of a
potential $G$, which is a fibrewise Legendre transform of a K\"ahler
potential, and is required to satisfy only open conditions---such `dual
potentials' appeared first in the toric case, in work of
Guillemin~\cite{G:kstv} and Abreu~\cite{Abreu}. We also describe the invariant
pluriharmonic functions and compute the Ricci form.

In subsections~\ref{s:rigid} and~\ref{s:ss} we introduce rigid and semisimple
hamiltonian torus actions respectively. In the case of circle actions (in
particular) these conditions originally arose as an Ansatz for the
construction of extremal K\"ahler metrics and K\"ahler--Einstein
metrics~\cite{calabi0,calabi1,hwang-singer,koisak,ped-poon,
christina,christina2,wang}: for semisimple rigid actions in momentum
coordinates, the Ricci form is linear in the matrix of inner products of the
Killing vector fields generating the action.  Our approach provides a natural
interpretation, particularly for the rigidity condition.

Subsection~\ref{s:ortho} is devoted to orthotoric K\"ahler metrics: in four
dimensions, these were introduced in~\cite{ACG0} and explicitly classified;
here we extend the definition and classification to all dimensions. The
decisive feature of orthotoric $2m$-manifolds is that they depend effectively
on $m$ functions of $1$ variable, rather than the $1$ function of $m$
variables (the dual potential $G$) that governs toric K\"ahler metrics in
general. This means that curvature conditions are (functional) ordinary
differential equations, rather than partial differential equations.

The central results of this paper can be found in section~\ref{s:chf}, where
we bring the work of sections~\ref{s:ham} and~\ref{s:hta} together. We first
prove that on a connected K\"ahler $2m$-manifold with hamiltonian $2$-form
$\phi$ and associated Killing vector fields $K_1,\ldots K_m$, there is an
integer $\ell$, with $0\leq\ell\leq m$, such that the span of $K_1,\ldots K_m$
is everywhere at most $\ell$-dimensional, but on a dense open set $K_1,\ldots
K_\ell$ are linearly independent. We show that $\ell$ roots of the {\it
momentum polynomial} $\Mp(t):=(-1)^m\pfaff(\phi - t\omega)$ are functionally
independent, the remainder being constant. (Here $\pfaff \psi =
\frac{1}{m!}{*(\psi^{\wedge m})}$ stands for the pfaffian of a $2$-form
$\psi$.) We call $\ell$ the {\it order} of $\phi$ and prove that K\"ahler
manifolds admitting a hamiltonian $2$-form of order $\ell$ are exactly those
admitting a local hamiltonian $\ell$-torus action such that
\begin{itemize}
\item the fibres are orthotoric;
\item the action is rigid;
\item the base is semisimple, with relative eigenvalues of a special form.
\end{itemize}

In section~\ref{s:ckmh} we study the curvature of our explicit metrics and
hence obtain classifications of extremal K\"ahler metrics with a hamiltonian
$2$-form, of weakly Bochner-flat K\"ahler metrics, and of Bochner-flat
K\"ahler metrics~\cite{bryant}.

To summarize, we have the following result.
\begin{result}
Let $(M,g,J,\omega)$ be a connected K\"ahler $2m$-manifold with a hamiltonian
$2$-form $\phi$ of order $\ell$. Then there are functions $\pF_1,\ldots
\pF_\ell$ of one variable such that on a dense open subset $M^0$ of $M$, the
K\"ahler structure may be written
\begin{align*}
g&=\sum_\xi\Mpn(\xi) g_{\xi}
+\sum_{j=1}^\ell \frac{\Mp'(\xi_j)}{\pF_j(\xi_j)} d\xi_j^2
+\sum_{j=1}^\ell \frac{\pF_j(\xi_j)}{\Mp'(\xi_j)}\Bigl(\sum_{r=1}^\ell
\sigma_{r-1}(\hat\xi_j)\theta_r\Bigr)^2,\\
\omega&=\sum_\xi \Mpn(\xi)
\omega_{\xi}
+\sum_{r=1}^\ell d\sigma_r\wedge \theta_r,\qquad\qquad\,
d\theta_r=\sum_\xi(-1)^r\xi^{\ell-r}\omega_{\xi},\\
J d \xi _j &= \frac{\pF_j (\xi _j)}{\Mp'(\xi_j)} \,
\sum _{r = 1} ^\ell \sigma _{r - 1} (\hat{\xi} _j) \,\theta_r,\qquad\qquad
\qquad J \theta_r =
(-1) ^r \,\sum_{j=1}^\ell \frac{\Mpc(\xi_j)}{\pF_j(\xi _j)} \xi_j^{\ell-r}
\,d \xi _j.
\end{align*}
Any K\"ahler metric of this form admits a hamiltonian $2$-form of
order $\ell$, namely
\begin{equation*}
\phi = \sum_\xi \xi\, \Mpn(\xi) \omega_\xi
+\sum_{r=1}^\ell (\sigma_r d\sigma_1 -d\sigma_{r+1})\wedge\theta_r.
\end{equation*}
In these expressions\textup:

\indent\llap{$\bullet$} $\sigma_r$ is the $r$th elementary symmetric
function of the non-constant roots $\xi_1,\ldots\xi_\ell$ of the
momentum polynomial $\Mp(t)$---so $\sigma_{\ell+1}=0$---and
$\sigma_{r-1} (\hat{\xi} _j)$ is the $(r-1)$st elementary symmetric
functions of the $\ell-1$ roots $\{\xi_k:k\neq j\}$\textup;

\indent\llap{$\bullet$} summation over $\xi$ denotes the sum over the
\emph{different} constant roots of the momentum polynomial and
$(g_\xi,\omega_\xi)$ is a positive or negative definite K\"ahler metric on a
manifold $S_\xi$ of the same \textup(real\textup) dimension $2m_\xi$ as the
$\xi$-eigenspace of $-J\circ\phi$\textup;

\indent\llap{$\bullet$} $\Mp(t)=\Mpn(t)\Mpc(t)$, where
$\Mpn(t)=\prod_{j=1}^\ell (t-\xi_j)$ and $\Mpc(t)=\prod_\xi
(t-\xi)^{m_\xi}$\textup; also $\Mp'(t)$ is the $t$-derivative of $\Mp(t)$, so
that $\Mp'(\xi_j)=\Mpc(\xi_j)\prod _{k \neq j} (\xi_j-\xi_k)$.

\smallbreak Now define polynomials $\check \Mpc(t)=\prod_\xi
(t-\xi)^{m_\xi-1}$, and $\hat \Mpc(t)=\prod_\xi (t-\xi)^{m_\xi+1}$.
Then we have the following special cases.

\begin{numlist}
\item $(g,J,\omega)$ is an extremal K\"ahler metric if
\begin{itemize}
\item for all $j$, $\pF_{j}''(t)=\check \Mpc(t)\bigl( \sum_{r=0}^{\check
m} a_r t^{\check m-r}\bigr)$, where $a_0,\ldots a_{\check m}$ are arbitrary
constants \textup(independent of $j$\textup) and $\check m=\ell+\sum_\xi 1 = m
- \sum_\xi (m_\xi-1)$\textup;
\item for all $\xi$, $\pm(g_\xi,\omega_\xi)$ has
$\Scal_{\pm g_\xi}=\mp\bigl(\sum_{r=0}^{\check m} a_r \xi^{\check m-r}\bigr)
/\prod_{\eta\neq\xi} (\xi-\eta)$.
\end{itemize}
The scalar curvature is constant if and only if $a_0=0$ and zero if and only
if also $a_1=0$.  An extremal K\"ahler metric with a hamiltonian $2$-form
arises in this way if the gradient of the scalar curvature is tangent to the
fibration defined by the $2$-form.

\item $(g,J,\omega)$ is weakly Bochner-flat if
\begin{itemize}
\item for all $j$, $\pF_j'(t)= \Mpc(t)\bigl(\sum_{r=-1}^{\ell} b_r
t^{\ell-r}\bigr)$, where $b_{-1},\ldots b_{\ell}$ are arbitrary constants
\textup(independent of $j$\textup)\textup;
\item for all $\xi$, $\pm(g_\xi,\omega_\xi)$ is K\"ahler--Einstein with
`K\"ahler--Einstein constant'
\begin{equation*}\notag
\frac1{m_\xi}\Scal_{\pm g_\xi}=\mp\sum_{r=-1}^\ell b_r\xi^{\ell-r}.
\end{equation*}
\end{itemize}
$(g,J,\omega)$ is K\"ahler--Einstein if and only if $b_{-1}=0$ and Ricci-flat
if and only if also $b_0=0$.  Any weakly Bochner-flat K\"ahler metric arises
in this way.

\item $(g,J,\omega)$ is Bochner-flat if
\begin{itemize}
\item for all $j$, $\pF_j(t)= \hat \Mpc(t)\bigl(\sum_{r=-2}^{\hat m} c_r
t^{\hat m-r}\bigr)$, where $c_{-2},\ldots c_{\hat m}$ are arbitrary constants
\textup(independent of $j$\textup) and $\hat
m=\ell-\sum_\xi 1=m-\sum_\xi (m_\xi+1)$\textup;
\item for all $\xi$, $\pm(g_\xi,\omega_\xi)$ has constant holomorphic
sectional curvature
\begin{equation*}\notag
\frac1{m_\xi(m_\xi+1)}\Scal_{\pm g_\xi}
=\mp\Bigl(\sum_{r=-2}^{\hat m} c_r\xi^{\hat m-r}\Bigr)
\prod_{\eta\neq\xi}(\xi-\eta).
\end{equation*}
\end{itemize}
$(g,J,\omega)$ has constant holomorphic sectional curvature if and only if
$c_{-2}=0$ and is flat if and only if also $c_{-1}=0$.  Any Bochner-flat
K\"ahler metric arises in this way.
\end{numlist}
\end{result}
This theorem follows from Theorems~\ref{t:one} and~\ref{t:two}, and
Propositions~\ref{p:ext},~\ref{p:wbf} and~\ref{p:bf} in the text below.  We
end by discussing hamiltonian $2$-forms of order $1$, and the classification
of conformally Einstein K\"ahler metrics~\cite{derd-masch}.  There are also
two appendices. In Appendix A, we relate hamiltonian $2$-forms to conformal
Killing forms, recently studied by Moroianu and Semmelmann~\cite{MS,sem}. In
Appendix B, we collect some Vandermonde identities, which we have used freely
in the paper.

We thank Uwe Semmelmann for discussing conformal Killing forms with us, and
Christina T\o nnesen-Friedman for her helpful comments and interest in this
work.

\tableofcontents

\vspace{-8mm}

\section{The curvature of a K{\"a}hler manifold}\label{s:ckm}

In this section we review some background material in order to fix notation,
and to present the notions of Bochner-flat, weakly Bochner-flat, and
conformally Einstein K\"ahler metrics. Our conventions mainly
follow~\cite{besse}.

\subsection{Riemannian curvature}

The curvature $R$ of a $n$-dimensional riemannian manifold $(M, g)$ is defined
by
\begin{equation*} R _{X, Y} Z = \nabla _{[X, Y]} Z - [\nabla _X, \nabla _Y] Z,
\end{equation*}
for all vector fields $X, Y, Z$, where $\nabla$ denotes the Levi-Civita
connection. It is a $2$-form with values in the adjoint bundle $AM$
(the bundle of skew endomorphisms of the tangent bundle $TM$) and satisfies the
{\it algebraic Bianchi identity}: $R _{X, Y} Z + R _{Y, Z} X + R _{Z, X} Y =
0$.  Via the metric $g$, $AM$ can be identified with the bundle $\Lambda ^2 M$
of $2$-forms and $R$ can be viewed as a section of $\Lambda ^2 M \otimes
\Lambda ^2 M$. Then, the algebraic Bianchi identity is equivalent to the
following two conditions:
\begin{numlist}
\item $R$ belongs the the symmetric part, $S^2 \Lambda ^2 M$, of $\Lambda ^2 M
\otimes \Lambda ^2 M$;

\item $R$ belongs to the kernel of the linear map, $\beta$, from $S^2\Lambda
^2 M$ to $\Lambda ^4 M$ determined by the wedge product.
\end{numlist}
${\mathcal R}M := \ker \beta\subseteq S^2\Lambda ^2 M$ is called the
bundle of (abstract) curvature tensors.

The {\it Ricci contraction} is the linear map $c$ from ${\mathcal R} M$ to the
bundle $S M$ of symmetric bilinear forms of $M$ sending $R$ to the
bilinear form $\Ric$ defined by $\Ric _{X, Y} = \trace (Z \to R_{X, Z} Y)$.
We thus obtain an orthogonal decomposition:
\begin{equation*} {\mathcal R} M = c ^* (S M) \oplus {\mathcal W} M,
\end{equation*}
where ${\mathcal W} M$, called the bundle of (abstract) Weyl tensors of $(M,
g)$, denotes the kernel of $c$ in ${\mathcal R} M$. Accordingly, the curvature
$R$ splits as $R = c ^* (h) + W$, where $W$ is the {\it Weyl tensor} of $(M,
g)$, whereas $h$ satisfies $c c ^* (h) = \Ric$.  For $n \geq 3$, $c$ is
surjective; its adjoint $c ^*$ is then injective and $h$ is determined by
\begin{equation*} h = \frac{\Scal}{2 n (n - 1)} g + \frac{\Ric _0}{n - 2},
\end{equation*}
where $\Scal$ is the {\it scalar curvature} of $g$, i.e., the trace of $\Ric$
with respect to $g$, and $\Ric _0$ denotes the traceless part of $\Ric$ (so
that $\Ric = \frac 1 n \Scal \,g + \Ric _0$), which is a section of $S_0 M$,
the bundle of symmetric traceless bilinear forms.  For $n=2$, $c^*$ has kernel
$S_0M$ so that $\Ric_0=0$ and the tracefree part of $h$ is undetermined.

Finally, the curvature $R$, viewed as a symmetric endomorphism of
$\Lambda^2M$ using $g$, splits into three orthogonal pieces as
follows:
\begin{equation} \label{decR}
R = \frac{\Scal}{n (n - 1)}\Id \restr{\Lambda ^2 M}
+ \frac{1}{n - 2} \{\Ric _0,\cdot\}  + W,
\end{equation}
where $\{\Ric _0,\cdot\}$ acts on $\psi\in\Lambda^2M$ to give the
anticommutator $\{ \Ric_0, \psi \} := \Ric_0 \circ\psi + \psi\circ \Ric_0$
of $\Ric_0$ and $\psi$, which are viewed, via $g$, as endomorphisms,
respectively symmetric and skew, of $TM$. When $n=2$, the second term is zero.

Each piece of (\ref{decR}) is an element of ${\mathcal R} M$. The
corresponding subbundle of ${\mathcal R} M$ is associated to an irreducible
representation of the orthogonal group $O (n)$, respectively the trivial
representation, the Cartan product ${\R} ^n \odot {\R} ^n$, the Cartan product
${\mathfrak o} (n) \odot {\mathfrak o} (n)$, where ${\mathfrak o} (n) \cong
\Lambda ^2 {\R} ^n$ denotes the Lie algebra of $O (n)$. (The Cartan product of
two irreducible representations with dominant weights $\lambda _1$ and
$\lambda _2$ is the irreducible sub-representation of the tensor product with
dominant weight $\lambda _1 + \lambda _2$.)

\subsection{The Bochner tensor of a K\"ahler manifold}

Let $(M, g,J,\omega)$ be a K{\"a}hler manifold of dimension $n=2m$.  By
definition, $J$ is an orthogonal complex structure which is parallel with
respect to the Levi-Civita connection $\nabla$. The K{\"a}hler form $\omega$
is defined by $\omega (X, Y) = g (JX, Y)$. The {\it Ricci form} $\rho$ and its
primitive part $\rho _0$ are defined in a similar way: $\rho (X, Y) = \Ric
(JX, Y)$ and $\rho _0 (X, Y) = \Ric _0 (JX, Y)$. The Ricci tensor is
$J$-invariant, and so $\rho$ and $\rho_0$ are $2$-forms.
We denote by
\begin{equation*}\notag
\Lambda ^2 M = \Lambda ^{J, +}M \oplus \Lambda ^{J, -} M,
\end{equation*}
the (orthogonal) decomposition of $\Lambda ^2 M$ into its $J$-invariant part,
$\Lambda ^{J, +}M$, and its $J$-anti-invariant part, $\Lambda ^{J, -} M$. The
riemannian curvature $R$ has values in $\Lambda ^{J, +}M$ and therefore acts
trivially on $\Lambda ^{J, -} M$. More generally, we call an element of
${\mathcal R} M$ {\it k{\"a}hlerian} if it acts trivially on $\Lambda ^{J, -}
M$. The set of abstract k{\"a}hlerian curvature tensors is a vector subbundle
of ${\mathcal R} M$, denoted by ${\mathcal K} M$; thus ${\mathcal K} M$ is the
kernel of the linear map from $S^2 \Lambda ^{J, +} M$ to $\Lambda ^4 M$
determined by the wedge product.

The curvature tensor $R$ of a K{\"a}hler manifold $(M, g, J)$ is a section of
${\mathcal K} M$, but in general none of its components in~\eqref{decR} is.
Indeed, the first component of $R$ in \eqref{decR} is only an element of
${\mathcal K} M$ if it is zero or $n=2$, while the second component is only an
element of $\mathcal K M$ if it is zero or $n=4$.  We define the {\it Bochner
tensor} $\WK$ to be the orthogonal projection of the third component, the Weyl
tensor $W$, onto ${\mathcal W} M \cap {\mathcal K} M$.  We thus obtain a new
decomposition of the curvature $R$ inside $\mathcal K M$.
\begin{equation} \label{decRK} \begin{split}
R &= \frac{\Scal}{2m(m+1)}(\Id\restr{\Lambda^{J,+}M} + \omega\otimes\omega) \\
& \quad + \frac{1}{m + 2} \left (\{\Ric_0,\cdot\}\restr{\Lambda^{J,+}M}
+ \rho _0 \otimes \omega + \omega \otimes \rho _0 \right ) \\
& \quad + \WK.
\end{split} \end{equation}
Here ${}\restr{\Lambda ^{J, +} M}$ has to be interpreted as the orthogonal
projection $\psi\mapsto\psi^{J,+}$ of $\Lambda ^2 M$ onto its $J$-invariant
part $\Lambda ^{J, +} M$, and $\rho_0\otimes\omega$ acts on
$\psi\in\Lambda^2M$ to give $\g{\rho_0,\psi}\omega$, where the inner product
on $2$-forms is normalized so that $\g{\omega,\omega} = m = n/2$.

The three pieces of $R$ appearing in (\ref{decRK}) are sections of subbundles
of ${\mathcal K} M$ associated to irreducible representations of the unitary
group $U (m)$, viewed as a subgroup of $O (n)$, namely: the trivial
representation, the Cartan product ${\C} ^m \odot {\C} ^m$, the Cartan product
$\mathfrak{su} (m) \odot \mathfrak{su} (m)$ respectively, $\mathfrak {su} (m)$
being the Lie algebra of $SU (m)$.

If $\Scal$ is a positive constant, then the first component of $R$ in
\eqref{decRK} agrees with the curvature of the complex projective space $\C
P^m$ with the Fubini--Study metric of holomorphic sectional curvature equal to
$\frac{\Scal}{m(m+1)}$.

The second component in~\eqref{decRK} agrees with the second component in
\eqref{decR} when $n=4$, since then $\{ \Ric_0, \psi^{J,-} \}=[J \psi^{J,-},
\rho _0] = 0$, whereas $\{ \Ric_0, \psi^{J,+} \}
=\g{\psi,\rho_0}\omega+\g{\psi,\omega}\rho_0$. The four dimensional case is
also special because the Weyl tensor $W$ splits into selfdual and antiselfdual
parts as $W = W ^+ + W ^-$, and on a K\"ahler $4$-manifold the selfdual part
is identified with the scalar curvature by
\begin{equation*}
W ^+ = \frac{\Scal}{12}\Big(\frac{3}{2}\omega\otimes\omega - \Id\restr{\Lambda
^{J, +} M}\Bigr).
\end{equation*}
Bringing together $W ^+$ and the scalar part of $R$ in \eqref{decR}, we deduce
that $\WK = W^-$.

In higher dimensions, the Weyl tensor $W$ of a K{\"a}hler manifold splits
into three pieces: one is the Bochner tensor $\WK$, while the other two are
identified with $\Ric_0$ and $\Scal$. In other words, on a K{\"a}hler manifold
of dimension $n\geq6$ the information given by the riemannian curvature is
already contained in the Weyl tensor; in particular, for $n \geq 6$, a locally
conformally flat K{\"a}hler metric is flat.

\subsection{The differential Bianchi identity in K{\"a}hler geometry}

The differential Bianchi identity
\begin{equation} \label{bianchi}
\nabla _X R _{Y, Z} + \nabla _Y R _{Z, X} + \nabla _Z R _{X, Y} = 0,
\end{equation}
easily implies the following one, known as the {\it Matsushima identity}:
\begin{equation} \label{mat} (\delta R) _{JX}  = - \nabla _X  \rho.
\end{equation}
(We specialize (\ref{bianchi}) by $X = e _j$, $Y = J e _j$, where $\{ e _j \}$
is a local, $J$-adapted, orthonormal frame and we observe that $\rho =
\frac{1}{2} \sum _{j = 1} ^n R _{e _j, J e _j}$; here we define $(\delta R)
_{JX}:=-\sum_{j=1}^n \nabla_{e_i}R_{e_i,JX}$.)  The Matsushima identity
immediately implies that the Ricci tensor of a K{\"a}hler manifold is parallel
if and only if the curvature is co-closed, as a $2$-form with values in
$\Lambda ^2 M$. The Ricci form $\rho$ may also be expressed as
$\rho(X,Y)=\frac{1}{2} \sum_{j=1} ^n \g{R _{X,Y} e_j,J e_j}$, and so it is
closed by~\eqref{bianchi}. Hence from \eqref{decRK} and \eqref{mat}, we infer
the following expression for the codifferential of the Bochner tensor:
\begin{equation} \label{MT} \begin{split}
(\delta \WK) _{JX} = - \frac{m}{m + 2} \, \nabla _X \rho _0
& - \frac{1}{2 (m+1) (m+2)} \, d \Scal (X) \, \omega \\
& + \frac{m}{4 (m+1) (m+2)} \,(d \Scal \wedge JX - d^c \Scal \wedge X).
\end{split} \end{equation}
Here $d^c= J\circ d$, $X$ is any vector field, and we identify vector fields
and $1$-forms via $g$.  In view of this identity, we introduce a {\it
normalized} Ricci form $\nrho$ defined by
\begin{equation*}
\nrho = \rho _0 + \frac{\Scal}{2 m (m + 1)} \, \omega.
\end{equation*}
Then, identity \eqref{MT} reduces to
\begin{equation} \label{MT1}
\frac{m+2}{m}(\delta \WK) _{JX} =
-\nabla _X \nrho + \frac{1}{2} (ds \wedge JX - d^c s \wedge X),
\end{equation}
where the {\it normalized scalar curvature} $s = \frac{\Scal}{2 (m + 1)}$ is
the trace of $\nrho$ with respect to $\omega$: $s = \g{\nrho, \omega}$.

\begin{defn}
A K{\"a}hler manifold $(M, g, J)$ is called {\it Bochner-flat} (or {\it
Bochner--K\"ahler}) if the Bochner tensor vanishes, $\WK=0$, and {\it weakly
Bochner-flat} if the Bochner tensor is co-closed, $\delta \WK = 0$.
\end{defn}
By \eqref{MT1}, a K{\"a}hler manifold is weakly Bochner-flat if and
only if it satisfies the following {\it weak Einstein condition}:
\begin{equation} \label{MT2}
\nabla _X \nrho = \frac{1}{2} (ds \wedge JX - d^c s \wedge X).
\end{equation}

\subsection{Conformally Einstein K\"ahler metrics}

A K\"ahler manifold $(M,g,J,\omega)$ of dimension $n=2m\geq 4$ is said to be
{\it conformally Einstein} if there is a nonvanishing function $\tau$ such
that $\tilde g:=\tau^{-2}g$ is an Einstein metric, i.e., $\Ric^{\tilde
g}_0=0$.  A straightforward and standard computation of the conformal change
of the Ricci tensor shows that $g$ is conformally Einstein with conformal
factor $\tau$ if and only if
\begin{equation} \label{confEin}
2(m-1)\nabla Jd\tau = - \tau\rho + \lambda \omega
\end{equation}
for some function $\lambda$---the trace of this equation then determines
that
\begin{equation*}
\lambda = -\frac{m-1}m \Delta \tau + \frac1{2m} \Scal\,\tau,
\end{equation*}
where $\Delta\tau = -\trace_g\nabla d\tau =-\g{dd^c\tau,\omega}$.

We recall that a hamiltonian vector field $K=J\grad_g f$ is Killing if and
only if it preserves $J$, if and only if the hessian $\nabla d f$ is
$J$-invariant, if and only if $\nabla Jdf=\frac12 dd^cf$, in which case $f$ is
said to be a {\it Killing potential}.

Clearly equation~\eqref{confEin} implies that $2(m-1)\nabla Jd\tau =
(m-1)dd^c\tau$, so that by differentiating~\eqref{confEin}, we obtain
$d\tau\wedge\rho-d\lambda\wedge\omega=0$ and hence $d\tau\wedge
d\lambda\wedge\omega=0$. We shall say that $g$ is {\it strongly}
conformally Einstein if $d\tau\wedge d\lambda=0$; this is automatic if
$n\geq 6$ since the wedge product with $\omega$ is then injective on
$2$-forms.

The fact that conformally Einstein K\"ahler metrics are strongly conformally
Einstein in $6$ or more dimensions was first observed by Derdzi\'nski and
Maschler~\cite{derd-masch}, who used this to obtain an explicit description of
such metrics. A key step is essentially equivalent to the following.

\begin{lemma}\textup{\cite{derd-masch}} Suppose that $g$ is strongly
conformally Einstein, with conformal factor $\tau$. Then on the
open set where $d\tau$ is nonzero,
\begin{equation}\label{eq:dm}
2\nabla Jd\tau = p \,\omega + q \,d\tau\wedge d^c\tau
\end{equation}
for some functions $p,q$ with $dp\wedge d\tau=0$.
\end{lemma}
\begin{proof}
On the open set where $d\tau$ is nonzero, we may write $d\lambda=\lambda_\tau
d\tau$, so that $d\tau\wedge(\rho-\lambda_\tau\omega)=0$.  It follows that the
$J$-invariant $2$-form $\rho-\lambda_\tau\omega$ is equal to $f d\tau\wedge
d^c\tau$ for some function $f$. Therefore:
\begin{equation*}
2(m-1)\nabla Jd\tau = (\lambda-\tau \lambda_\tau)\omega - \tau f \,d\tau\wedge
d^c\tau
\end{equation*}
and clearly $d(\lambda-\tau \lambda_\tau)\wedge d\tau=0$.
\end{proof}
\begin{rem} More generally, the conclusions of this lemma hold on the
open set where $d\tau, \xi_\tau\neq 0$ if $2\nabla Jd\tau =
\xi(\tau)\rho+\eta(\tau)\omega$, with essentially the same
proof~\cite{derd-masch}.
\end{rem}

In order to interpret the work of Derdzi\'nski and Maschler in the present
work, we reformulate equation~\eqref{eq:dm}. We first note that if $\tau$ is
any function satisfying~\eqref{eq:dm}, for some functions $p,q$ with $dp\wedge
d\tau=0$, then in fact we have
\begin{equation*}
d(|d\tau|^2)\wedge d\tau=0,\qquad dq\wedge d\tau=0,\qquad\text{and}
\qquad p=\frac{a}{a\tau+b}|d\tau|^2
\end{equation*}
for some constants $a$ and $b$ not both zero. Indeed,
contracting~\eqref{eq:dm} with $Jd\tau$ we obtain
$d(|d\tau|^2)=(p+q|d\tau|^2)d\tau$ which gives the first two observations.
Hence $dd^c\tau= f|d\tau|^2\omega+q\, d\tau\wedge d^c\tau$, where $df\wedge
d\tau=0$. The exterior derivative of this equation gives
$|d\tau|^2(df+f^2d\tau)\wedge\omega=0$, so that $f=a/(a\tau+b)$.

\begin{lemma} A Killing potential $\tau$ satisfies the equation
\begin{equation*}
dd^c\tau=\frac{a}{a\tau+b}|d\tau|^2\omega+q\, d\tau\wedge d^c\tau
\end{equation*}
\textup(for some function $q$ and constants $a,b$ not both zero\textup) if
and only
if the $2$-form $\chi:=(a\tau+b)d\tau\wedge d^c\tau/|d\tau|^2$ satisfies
\begin{equation}\label{DM}
\nabla_X\chi = \frac a2 (d\tau\wedge JX - d^c\tau\wedge X).
\end{equation}
\end{lemma}
\begin{proof}
If $\chi=(a\tau+b)d\tau\wedge d^c\tau/|d\tau|^2$, then
\begin{multline*}
\nabla_X\chi = a\,d\tau(X)\frac{d\tau\wedge d^c\tau}{|d\tau|^2}
-(a\tau+b)\g{\iota_X dd^c\tau,d^c\tau}\frac{d\tau\wedge d^c\tau}{|d\tau|^4} \\
+(a\tau+b)\frac{\iota_{JX}dd^c\tau\wedge d^c\tau+d\tau\wedge\iota_X dd^c\tau}
{2|d\tau|^2}.
\end{multline*}
This can only equal $\frac a2 (d\tau\wedge JX - d^c\tau\wedge X)$, if
$dd^c\tau$ is of the form $f|d\tau|^2\omega+q\,d\tau\wedge d^c\tau$ for some
functions $f,q$, in which case we obtain
\begin{equation*}
\nabla_X\chi=(a-(a\tau+b)f)d\tau(X)\frac{d\tau\wedge d^c\tau}{|d\tau|^2}
+\frac12 (a\tau+b)f\,(d\tau\wedge JX - d^c\tau\wedge X).
\end{equation*}
The result is now immediate.
\end{proof}

\section{Hamiltonian $2$-forms}\label{s:ham}

In this section we introduce the notion of a hamiltonian $2$-form and develop
the most basic general properties and the simplest examples.

\subsection{Hamiltonian $2$-forms and Killing vector fields}

The definition of hamiltonian $2$-forms is motivated both by weakly
Bochner-flat K\"ahler manifolds and by strongly conformally Einstein K\"ahler
manifolds. The reason for the terminology will shortly become apparent.

\begin{defn} Let $\phi$ be any (real) $J$-invariant $2$-form on the K{\"a}hler
manifold $(M, g, J,\omega)$.  We say $\phi$ is {\it hamiltonian} if there
is a function $\sigma$ on $M$ such that
\begin{equation} \label{ham}
\nabla _X \phi = \frac12 (d \sigma \wedge JX - d ^c \sigma \wedge X)
\end{equation}
for any vector field $X$. When $M$ is a Riemann surface, we require in
addition that $\sigma$ is a Killing potential.
\end{defn}
It follows immediately from the definition that $d\sigma=d\trace\phi$, where
$\trace \phi = \g{\phi, \omega}$ is the trace of $\phi$ with respect to
$\omega$, so without loss of generality, we may take $\sigma=\trace\phi$.  The
defining equation for hamiltonian $2$-forms is therefore linear.  Note that,
for a general hamiltonian $2$-form $\phi$, $A=\phi+\sigma\omega$ is closed,
$dA=0$.

\begin{ex}\label{ex1}
On any K\"ahler manifold, any $J$-invariant parallel $2$-form is hamiltonian.
In particular, a constant multiple of the K\"ahler form $\omega$ is
hamiltonian. It follows that if $\phi$ is hamiltonian, then so is
$\phi_t:=\phi-t\omega$ for any constant $t$.
\end{ex}

\begin{ex}
We shall be particularly interested in the hamiltonian $2$-forms arising from
the following immediate consequence of equation~\eqref{MT2}.
\begin{prop} \label{wbfcrux}
A K{\"a}hler manifold of dimension $2m\geq 4$ is weakly Bochner-flat if and
only if the normalized Ricci form $\nrho = \rho _0 + \frac{1}{m} s\,
\omega=\rho-s\,\omega$ is hamiltonian.
\end{prop}
\end{ex}

\begin{ex} In view of equation~\eqref{DM}, we also have the following result.
\begin{prop}\label{scecrux}
On any strongly conformally Einstein K\"ahler manifold of dimension $2m\geq
4$, with conformal factor $\tau$, there are constants $a,b$ not both zero,
such that $\chi=(a\tau+b)d\tau\wedge Jd\tau/|d\tau|^2$ is hamiltonian on the
open set where $d\tau$ is nonzero.
\end{prop}
\end{ex}

The equation for hamiltonian $2$-forms is overdetermined.  By
differentiating and skew-symmetrizing \eqref{ham}, we get
\begin{equation} \label{maincurv} \begin{split}
R _{X, Y}\cdot\phi = [R _{X, Y}, \phi] &=
\frac 1 2 (\nabla _Y d\sigma \wedge JX - \nabla _X d\sigma \wedge JY \\
& \quad\; - J\nabla _Y d\sigma \wedge X + J\nabla _X d\sigma \wedge Y).
\end{split} \end{equation}
This formula underlies most of the basic theory of hamiltonian $2$-forms.  In
particular, we shall use it to explain the use of the term ``hamiltonian'' in
this context.

To do that, we first recall that the {\it pfaffian} of a $2$-form
$\phi$ is defined by
\begin{equation} \label{pfaff}
\pfaff \phi=\frac1{m!}{*(\phi\wedge\cdots\wedge\phi)},
\end{equation}
where $*$ denotes the Hodge star operator. The normalization is chosen so that
$\pfaff\omega=1$ and thus $\phi \wedge \cdots \wedge \phi =
(\pfaff\phi)\omega\wedge\cdots\wedge\omega$. We let $\phi_t=\phi-t\omega$ as
in Example~\ref{ex1} above, and (following Bryant~\cite{bryant}) define the
{\it momentum polynomial} of $\phi$ to be
\begin{equation}\label{mompoly}
\Mp(t):=(-1)^m\pfaff \phi_t
= t^m - (\trace\phi) \,t^{m-1} + \cdots + (- 1) ^m \pfaff\phi.
\end{equation}

\begin{prop} \label{prophamilton}
If $\phi$ is a hamiltonian $2$-form, then the functions $\Mp(t)$ on $M$
\textup(for each $t\in\R$\textup) are Poisson-commuting hamiltonians for
Killing vector fields $K(t) := J \grad_g \Mp(t)$ which preserve $\phi$. In
particular, the vector fields $K(t)$ all commute.
\end{prop}
\begin{proof} We first prove that $K := J \grad_g \sigma$ is Killing,
i.e., $\nabla d\sigma$ is $J$-invariant.  Since $R_{X,Y}$ is $J$-invariant in
$X$ and $Y$, equation~\eqref{maincurv} implies that
\begin{gather} \label{eq:S=0}
S (X) \wedge JY - J S (X) \wedge Y - S(Y) \wedge JX + J S (Y) \wedge X = 0,\\
\tag*{where}
S (X) = \nabla _X d\sigma + J \nabla _{JX} d\sigma.
\end{gather}
Contracting~\eqref{eq:S=0} with a vector field $Z$ and taking the trace
over $Y$ and $Z$ yields $2(1-m) J S(X)=0$ and hence $\nabla d\sigma$ is
$J$-invariant---by definition when $m=1$.

We now show that the other hamiltonian vector fields are Killing. To do this
we differentiate $\pfaff\phi_t$, using the fact that $\phi_t$ is hamiltonian
with $\trace\phi_t=\trace\phi-mt$ and hence $d \trace\phi_t=d\sigma$.
Therefore, from \eqref{ham} and \eqref{pfaff}, we get
\begin{equation}\label{dpfphi}
d \pfaff\phi_t = \frac{1}{(m-1) !}
{* (J d \sigma \wedge \phi_t \wedge \cdots \wedge \phi_t)}.
\end{equation}
Using \eqref{ham} again, we then obtain
\begin{equation} \label{nabladpf}\begin{split}
\nabla_X d \pfaff\phi_t & =  \frac{1}{(m - 1) !}
{* (\nabla _X J d \sigma \wedge \phi_t \wedge \cdots \wedge \phi_t)} \\
& \quad + \frac 1{2\,(m-2) !} {* ( JX \wedge d \sigma \wedge J d \sigma
\wedge \phi_t \wedge \cdots \wedge \phi_t)}.
\end{split} \end{equation}
The second term on the right hand side is automatically $J$-invariant, while
the first one is also $J$-invariant since $K$ is Killing. Hence $J\grad_g
\Mp(t)$ is Killing for all $t$.

It remains to prove that the Killing vector fields preserve $\phi$ and that
their momentum maps Poisson-commute. Contracting equation~\eqref{dpfphi} with
$Jd\sigma$, we deduce that $\g{Jd \pfaff\phi_t,d\sigma}=0$ and hence $K(t)$
preserves $\sigma$ for all $t$. It follows that $\mathcal L_{K(t)}\phi=\mathcal
L_{K(t)}(\phi_t+\sigma\omega)= d\,\iota_{K(t)}(\phi_t+\sigma\omega)$, since
$\phi_t+\sigma\omega$ is closed.

Now equation~\eqref{dpfphi} also implies that
\begin{equation*}
\phi_t\bigl(\grad_g (\pfaff\phi_t),\cdot\bigr)=\frac1{m!}\sum_{j=1}^{\smash
2m} \g{Jd\sigma,{*\iota_{e_j}(\phi_t\wedge\cdots\wedge\phi_t)}}\eps_j
=(\pfaff\phi_t) Jd\sigma
\end{equation*}
(using a local frame $e_j$ with dual frame $\eps_j$) so that
\begin{equation}\label{Krelation}
\phi_t(J\grad_g \Mp(t),\cdot)=-\Mp(t) d\sigma.
\end{equation}
Hence $\iota_{K(t)}(\phi_t+\sigma\omega)=-d(\sigma\,\Mp(t))$ is closed and so
$K(t)$ preserves $\phi$.  It follows that $K(t)$ preserves $\Mp(s)$ for all
$s,t\in\R$ and
$\{\Mp(s),\Mp(t)\}=\g{J\,K(s),K(t)}=d(\Mp(s))(K(t))=0$ for
all $s,t$; thus $\Mp(s)$ and $\Mp(t)$ Poisson-commute.
\end{proof}
Obviously $\Mp(t)$ is a Killing potential for all $t$ if and only if its
coefficients are all Killing potentials.

\subsection{The connection on $2$-jets of hamiltonian $2$-forms}

We have noted already that the equation for hamiltonian $2$-forms is
overdetermined. In fact it has finite type, i.e., the space of local
solutions is finite dimensional, being given by the parallel sections
for a connection on the $2$-jets of hamiltonian $2$-forms.

\begin{prop} If $\phi$ is a hamiltonian $2$-form then
\begin{gather}
\label{nablaphi}
\nabla \phi + {\frac12} (K \wedge \Id + JK \wedge J) = 0\\
\label{nablaK}
\nabla K + {\frac{\smash1}{2m}}
\bigl(2u\,\omega -  J\{\rho,\phi\} - 2R(\phi) )=0\\
\label{du}
du+\rho(K)=0.
\end{gather}
\textup(Here, as elsewhere, we identify vectors with $1$-forms and bilinear
forms with endomorphisms using $g$, and we recall that $\{\cdot,\cdot\}$
denotes the anticommutator.\textup)

Thus $(\phi,K,u)$---with $K=J\grad_{\smash g}\sigma$, $u=\frac12\Delta\sigma$
and $\sigma=\trace\phi$---is parallel with respect to a natural covariant
derivative $\mathcal D$ on $\Lambda^{J,+}M\oplus TM\oplus M{\times}\R$.

The integrability condition $F^{\mathcal D}\cdot(\phi,K,u) = 0$ is equivalent
to the equations
\begin{gather}
\label{eq:curv1}
m[R(\psi),\phi]-[R(\phi)+\tfrac12J\{\rho,\phi\},\psi^{J,+}]=0\\
\label{eq:curv2}
(m+1) R_{K,X}-\rho(K,X)\omega
+\tfrac12J\{\rho,K\wedge X^{J,+}\}-\tfrac12J\{\nabla_X\rho,\phi\}
-\nabla_XR(\phi)=0\\
\label{eq:curv3}
-m\nabla_K\rho + [R(\phi), \rho]= 0.
\end{gather}
for any $2$-form $\psi$ and vector field $X$. Note that~\eqref{eq:curv1}
with $\psi=\omega$ gives $[\rho,\phi]=0$.
\end{prop}
\begin{proof}
Equation~\eqref{nablaphi} is immediate by definition.
Contracting~\eqref{maincurv} with a vector field $Z$ and taking the trace over
$Y$ and $Z$ gives ${\textstyle\sum_{j=1}^{2m}
R_{X,e_j}\phi(e_j)}+\phi(\Ric(X)) =-\tfrac12(\Delta\sigma)JX-m
J\nabla_Xd\sigma$ (for a local orthonormal frame $e_j$) and~\eqref{nablaK} is
the $J$-invariant part of this (as $\nabla d\sigma$ is $J$-invariant). Since
$K$ is a Killing vector field, $\nabla_X \nabla K = R_{K,X}$~\cite{kostant},
and~\eqref{du} is obtained by contracting this with $\omega$.

The first integrability condition~\eqref{eq:curv1} follows
from~\eqref{maincurv} by substituting for $\nabla K=J\nabla d\sigma$.
Differentiating~\eqref{nablaK} using $\nabla_X\nabla K = R_{K,X}$ and the
equations for $\nabla_X\phi$ and $du(X)$ gives~\eqref{eq:curv2}.  Finally,
from equation~\eqref{du}, $0=d(\rho(K))=\mathcal L_{K}\rho=\nabla_K\rho -
[\nabla K, \rho]$, which yields~\eqref{eq:curv3} by substituting for $\nabla
K$.

The three components of $F^{\mathcal D}\cdot(\phi,K,u)$ are the left hand
sides of~\eqref{eq:curv1}--\eqref{eq:curv3} divided by $m$, after applying the
isomorphism ${\rm alt}\colon T^*M\otimes \Lambda^2M \to \Lambda^2M \otimes TM$
to~\eqref{eq:curv2}.
\end{proof}
\begin{rem}
It follows that hamiltonian $2$-forms enjoy the properties of parallel
sections, such as unique continuation, extendibility to submanifolds of
codimension at least two, and an upper bound, here equal to $m^2+2m+1$, on the
dimension of space of hamiltonian $2$-forms.
\end{rem}

We now expand the curvature $R$ using~\eqref{decRK}, which may be rewritten
\begin{equation}\label{eq:newdecRK}
R(\psi) = \WK(\psi) -J\{\hat\rho,\psi^{J,+}\}
+\g{\hat\rho,\psi}\omega+(\trace\psi)\hat\rho
\end{equation}
for any $2$-form $\psi$, where $\hat\rho = \frac1{m+2}\bigl(\nrho+\tfrac12
s\,\omega\bigr)= \frac1{m+2}\bigl(\rho-\tfrac12
s\,\omega\bigr)$. Then~\eqref{nablaK} reads
\begin{equation} \label{nablaK1} \begin{split}
\nabla K&=\frac12 J\{\hat\rho,\phi\}+\frac1{m}
\bigl(\WK(\phi) +(\trace\phi)\hat\rho-(\trace\hat\rho)\phi+
(\g{\hat\rho,\phi}-u)\omega\bigr).
\end{split}\end{equation}
Using $[\rho,\phi]=0$, equation~\eqref{eq:curv1} implies
\begin{equation} \label{curvphi}
[\WK(\psi),\phi]=\frac1m [\WK(\phi),\psi^{J,+}]
+J\bigl(\hat\rho_0\circ\psi^{J,+}\circ\phi_0-
\phi_0\circ\psi^{J,+}\circ\hat\rho_0\bigr)
\end{equation}
$(\hat\rho_0=\frac{\rho_0}{m+2})$. Equation~\eqref{eq:curv2} is complicated
when fully expanded. Instead we use the fact that the $J$-invariant part of
$K\wedge X$ is $-\nabla_X\phi$ to obtain
\begin{equation*}
R_{K,X}=\WK_{K,X} + J\{\hat\rho,\nabla_X\phi\}
-\g{\hat\rho,\nabla_X\phi}\omega-d\sigma(X)\hat\rho
\end{equation*}
and $du(X)=-\rho(K,X)=(m+2)\g{\hat\rho,\nabla_X\phi}
+d\sigma(X)\trace\hat\rho$.  Substituting these into the covariant derivative
of~\eqref{nablaK1} (using $\nabla_X\nabla K = R_{K,X}$ as before) we have
\begin{equation}\label{eq:WKK1}
\begin{split}
\WK _{K, X} - \frac{1}{m} \nabla _X (\WK (\phi)) &=
\frac12J\{\nabla_X\hat\rho_0,\phi_0\}
+\frac1m\g{\nabla_X\hat\rho_0,\phi_0}\omega-\frac1m ds(X)\phi_0\\
&-\frac12J\{\hat\rho_0,\nabla_X\phi_0\}
-\frac1m\g{\hat\rho_0,\nabla_X\phi_0}\omega+\frac1m d\sigma(X)\rho_0\
\end{split} \end{equation}
The important point we shall need later is that the right hand side vanishes
if $\nrho$ is a constant linear combination of $\phi$ and $\omega$.

\subsection{The differential system in the weakly Bochner-flat case}
\label{s:dswbf}

On a weakly Bochner-flat K\"ahler manifold, the normalized Ricci form $\nrho$
is hamiltonian. We also want to study hamiltonian $2$-forms on
K\"ahler--Einstein manifolds. These cases can be considered together by
supposing that $(g,J,\omega)$ is a weakly Bochner-flat K\"ahler metric with a
hamiltonian $2$-form $\phi$ such that $\nrho$ is a constant linear combination
of $\phi$ and $\omega$.  We set $\nrho=(m+2)a\phi+b\omega$, and find
that~\eqref{nablaK1} may be written
\begin{equation*}\notag
\nabla K = \frac1m \WK(\phi) + a \,J\circ\phi^2 - (a\sigma+b)\phi+
\frac1m\Bigl(a(\sigma^2+\g{\phi,\phi})+b\sigma-\frac12\Delta\sigma\Bigr)\omega.
\end{equation*}
Let us put $\tau_0=-2a$, $\tau_1=-2(a\sigma+b)$ and
$\tau_2=\frac2m\bigl(a(\sigma^2+\g{\phi,\phi})+b\sigma-\frac12
\Delta\sigma\bigr)$. Then, as $d \g{\phi,\phi} = -2 \phi(K)$, we obtain the
following formulation of the system~\eqref{nablaphi}--\eqref{du}:
\begin{equation}\label{nablaphi2}\begin{split}
\nabla \phi &= -\frac12 (K \wedge \Id + JK \wedge J),\\
\nabla K &= \frac1m \WK(\phi),
-\frac12 \tau_0\,J\circ\phi^2 + \frac12\tau_1\,\phi + \frac12\tau_2\,\omega,\\
d\tau_2 &= -\tau_0 \,\phi(K) - \tau_1\, JK,\qquad
d\tau_1 = -\tau_0 \,JK,\qquad
d\tau_0 = 0.
\end{split}\end{equation}
When $\WK(\phi)=0$ this system yields an invariant polynomial---in particular
for $\tau_0\neq0$ and $\WK=0$, such a polynomial was found by
Bryant~\cite{bryant} and is the basis for his classification of Bochner-flat
K\"ahler metrics.
\begin{prop}\label{p:charpoly}
Let $(\phi,K,\tau_0,\tau_1,\tau_2)$ be a solution of~\eqref{nablaphi2} with
$\WK(\phi)=0$ and define a polynomial
\begin{equation}\label{eq:charpoly}
\CP(t):=(\tau_0 t^2+\tau_1 t+\tau_2)\Mp(t)-\g{K,K(t)}.
\end{equation}
Then $\CP(t)$ has constant coefficients. \textup(Recall that $K(t)=J\grad_g
\Mp(t)$, where $\Mp(t)=(-1)^m\pfaff\phi_t$ and $\phi_t=\phi-t\omega$.\textup)
\end{prop}
\begin{proof} Equation~\eqref{dpfphi} implies that
$\g{K,JX}\Mp(t) = \g{\phi_t(X),K(t)}$. Hence differentiating $(\tau_0
t^2+\tau_1 t+\tau_2)\Mp(t)$ along a vector field $X$, using the
system~\eqref{nablaphi2}, gives
\begin{align*}
&\g{K,\tau_0t\,JX}\Mp(t)+\g{K,\tau_0\,\phi(X)+\tau_1\, JX}\Mp(t)
+ (\tau_0 t^2+\tau_1 t+\tau_2)\g{JX,K(t)}\\
&\qquad=\g{(\tau_0t\,\phi_t-\tau_0\,\phi_t\circ\phi\circ J+\tau_1\,\phi_t
+\tau_0 t^2\,J+\tau_1 t\,J+\tau_2\,J)(X),K(t)}\\
&\qquad=\g{(-\tau_0\,J\circ\phi^2+\tau_1\,\phi+\tau_2\,J)(X),K(t)}
= 2\g{\nabla_X K, K(t)}.
\end{align*}
Now~\eqref{nabladpf} gives $\G{K,\nabla_X \bigl(K(t)\bigr)} =
\g{\nabla_X K, K(t)}$, which proves the proposition.
\end{proof}
Following Bryant, we refer to $\CP$ as the {\it characteristic polynomial}
of $(g,J,\omega,\phi)$.

\subsection{Complex projective, hyperbolic and euclidean space}\label{s:chsc}

A K\"ahler metric $g$ has constant holomorphic sectional curvature if and only
if it is Bochner-flat and K\"ahler--Einstein (and we require constant scalar
curvature when $m=1$). It follows from~\eqref{eq:curv1}--\eqref{eq:curv3} that
the connection $\mathcal D$ is flat in this case; hence on any simply
connected domain, the space of hamiltonian $2$-forms has dimension $(m+1)^2$.
Conversely, if $\mathcal D$ is flat, then $\rho_0=0$ (as $[\rho_0,\phi]=0$ for
all $\phi\in\Lambda^{J,+}M$); now~\eqref{curvphi} implies that $\WK=0$ (since
$[\WK(\phi),\psi]=0$ for all $\phi,\psi\in\Lambda^{J,+}M$); finally
\eqref{eq:WKK1} gives $ds=0$ (even if $m=1$), so $g$ has constant holomorphic
sectional curvature.

Hamiltonian $2$-forms on constant holomorphic sectional curvature manifolds
correspond to solutions of~\eqref{nablaphi2} with $\WK=0$, $\tau_0=0$ and
$\tau_1=-2s/m$. We first consider the case that $s$ is nonzero, i.e., up to
scale, the K\"ahler metric is the Fubini--Study metric of complex projective
space, or the Bergman metric of complex hyperbolic space. If we put
$\tau_2=2s\tau/m$, the system~\eqref{nablaphi2} becomes
\begin{equation}\label{nablaphi3}\begin{split}
\nabla \phi &= -\frac12 (K \wedge \Id + JK \wedge J)\\
\nabla K &= -\frac sm(\phi-\tau\omega)\\
d\tau&=JK.
\end{split}\end{equation}
The last two equations show that $\tau$ is a Killing potential for $-K$, and
that the hamiltonian $2$-form $\phi$ is completely determined by
$\tau$. Furthermore, the Kostant identity $\nabla_X\nabla K= R_{K,X}$ shows
that any Killing potential defines a hamiltonian $2$-form in this way.  Hence
there is a bijection $\phi\mapsto\frac1m \sigma-\frac 1{2s} \Delta\sigma$
(with inverse $\tau\mapsto \frac m {2s} dd^c\tau+\tau\omega$) from the space
hamiltonian $2$-forms to the space of Killing potentials, which may be
identified with the unitary Lie algebra $\mathfrak u(m+1)$ or $\mathfrak
u(m,1)$, using the Poisson bracket. We remark, though we shall not make use of
this, that the Lie bracket, Killing form and (monic, degree $m+1$)
characteristic polynomial can be computed and turn out to be given by
\begin{align*}
[(\phi,K,\tau),(\hat\phi,\hat K,\hat\tau)]&=\bigl(
[\phi,\hat\phi]+\tfrac m{s}{K{\wedge}\hat K}{}^{J,+}\!,
\iota_{\hat K}(\phi-\tau\omega)-\iota_K(\hat\phi-\hat\tau\omega),
\tfrac m{s}\omega(K,\hat K) \bigr)\\
\G{(\phi,K,\tau),(\hat\phi,\hat K,\hat\tau)}
&= \g{\phi,\hat\phi}+\tau\hat\tau+\tfrac m s \g{K,\hat K}\\
bP_{(\phi,K,\tau)}(t)&=-\tfrac m{2s}\CP(t)=
(t-\tau)\Mp(t)+\tfrac m{2s}\g{K,K(t)}.
\end{align*}

On a flat K\"ahler manifold (e.g., on complex euclidean space $\C^m$)
the system~\eqref{nablaphi2} reduces to
\begin{equation}\label{nablaphi4}\begin{split}
\nabla \phi &= -\frac12 (K \wedge \Id + JK \wedge J)\\
\nabla K &= \frac12\kappa\omega.
\end{split}\end{equation}
with $\kappa$ constant. Thus inside the space of hamiltonian $2$-forms we have
the parallel $2$-forms ($K=0$); modulo such parallel $2$-forms, we then have
the parallel vector fields ($\kappa=0$); then finally, the space of
hamiltonian $2$-forms on $\C^m$, modulo those with $K$ parallel, is one
dimensional, a representative element being $d t \wedge d ^c t$, where $t$ is
the distance squared to the origin. The characteristic polynomial is
now
\begin{equation*}
\CP(t)=\kappa\Mp(t)+\tfrac12\g{K,K(t)}
\end{equation*}
which has degree $m$ if $\kappa\neq 0$, degree $m-1$ if $\kappa=0$ and $K\neq
0$, and is zero if $K=0$.

We shall obtain an explicit description of the hamiltonian 2-forms on complex
projective, hyperbolic and euclidean space, with a given characteristic
polynomial $\CP(t)$, in section~\ref{s:bf} below.

\section{Hamiltonian torus actions}\label{s:hta}

We have seen that on a K\"ahler $2m$-manifold with a hamiltonian $2$-form,
there is a family of Poisson-commuting hamiltonian Killing vector fields
$K(t)=J\grad_g \Mp(t)$. Since $\Mp(t)$ is a monic polynomial of degree $m$,
the span of the $K(t)$ is at most $m$-dimensional. If they are not all zero,
then on an open set where the span has rank $\ell$, $1\leq \ell\leq m$, the
$K(t)$ generate a local action of an $\ell$-dimensional torus.

We next study hamiltonian $\ell$-torus actions in general. Our discussion is
independent of the theory of hamiltonian $2$-forms, but is strongly motivated
by it.  Roughly speaking, there are three aspects to the description of such
torus actions: first, the toric geometry of the fibres of the complexified
action; second, the geometry of the base of this action, the local K\"ahler
quotient; third, the way the fibre and base geometries fit together. In full
generality, these structures are quite difficult to handle. However, there is
a class of toric manifolds, called {\it orthotoric}, of K\"ahler quotients,
called {\it semisimple}, and of fibrations, called {\it rigid}, which are more
amenable to computation. It will turn out that the hamiltonian torus actions
induced by hamiltonian $2$-forms are always rigid with semisimple base and
orthotoric fibres.

\subsection{The Pedersen--Poon construction}

\begin{defn}\label{hamtorus}
A {\it local \textup(isometric\textup) hamiltonian $\ell$-torus action} on a
K{\"a}hler $2m$-manifold $(M,g,J,\omega)$ is an $\ell$-dimensional family of
holomorphic Killing vector fields $\boldsymbol K \in C^\infty(M,TM)\otimes
\R^{\ell*}$ which are linearly independent on a dense open set $M^0$ and
isotropic in the sense that $\omega(\boldsymbol K,\boldsymbol K)=0$.  The last
condition means that every component of $J\boldsymbol K$ is orthogonal to
every component of $\boldsymbol K$. It follows that $\ell\leq m$---if equality
holds, we say that $(M,g,J,\omega)$ is a {\it toric} K\"ahler manifold.
\end{defn}

For clarity, we write $K_r=\boldsymbol K(e_r)$ ($r=1,\ldots \ell$) for the
components of $\boldsymbol K$ with respect to a basis $e_r$ of
$\R^\ell$---this could also be interpreted as an abstract index notation. In
this subsection and the next two (only) we adopt the summation convention,
i.e., repeated indices imply contraction.

Since $\cL_{K_r}\omega=0$ for all $r$ and $\omega$ is closed, we have
$d(\iota_{K_r}\omega)=0$ and $\iota_{[K_r,K_s]}\omega=-d(\omega(K_r,K_s))=0$.
Furthermore, since $\cL_{K_r}J=0$ for all $r$ and $J$ is integrable, we have
$\cL_{JK_r}J=0$ and $[JK_r,JK_s]=J[JK_r,K_s]=0$.

\begin{rem} Definition~\ref{hamtorus} can be extended to almost hermitian
manifolds, but if $\omega$ is not closed, we assume a priori that
$[K_r,K_s]=0$ for all $r,s$, while if $J$ is not integrable, we assume that
$[JK_r,JK_s]$ is in the span of $J\boldsymbol K$ for all $r,s$.
\end{rem}

To obtain a local description of these metrics, valid near any point in $M^0$,
we may assume that $\boldsymbol K$ generates a free $\ell$-torus action, so
that $M$ is a principal $\ell$-torus bundle over a $(2m-\ell)$-dimensional
manifold $B$, and that the foliation generated by $\boldsymbol K, J\boldsymbol
K$ descends to a fibration of $B$ over a $2(m-\ell)$-dimensional manifold $S$.

Since $J$ is integrable and $\boldsymbol K$-invariant, the components of
$J\boldsymbol K$ are holomorphic vector fields, so that $S$ is a complex
manifold.

Further, since $\omega$ is closed and $\boldsymbol K$-invariant, we may
locally write $\iota_{\boldsymbol K}\omega=-d\boldsymbol\sigma$ where
$\boldsymbol\sigma\colon M\to\R^{\ell*}$ is a $\boldsymbol K$-invariant {\it
momentum map} for the torus action. Thus we may locally
identify $B$ with $S\times U$, where $U$ is an open subset of $\R^{\ell*}$ and
$\boldsymbol\sigma$ is given by projection to $U$. $S$ is then the {\it
K\"ahler quotient} of $M$: it is a complex manifold equipped with a
family of compatible K\"ahler structures parameterized by $U$.

It is useful to split the exterior derivative on $B$ into horizontal and
vertical parts:
\begin{equation*}\notag
d\alpha=d_h\alpha+d\sigma_r\wedge\cL_{\partial/\partial\sigma_r}\alpha.
\end{equation*}
(Note $\cL_{\partial/\partial\sigma_r}$ commutes with $d_h$. We write
$\partial\alpha/\partial\sigma_r$ as a shorthand for
$\cL_{\partial/\partial\sigma_r}\alpha$.)

Let $\boldsymbol\theta\in C^\infty(M,\Lambda^1M)\otimes\R^\ell$ be the
connection $1$-forms dual to $\boldsymbol K$ and let $J\boldsymbol
\theta=-\boldsymbol\theta\circ J$; thus $\theta_r(K_s)=\delta_{rs}$ and
$\theta_r$ vanishes on the horizontal distribution of $M\to B$ .  (We may
locally write $\theta_r=dt_r+\alpha_r$ where $\boldsymbol t\colon M\to\R^\ell$
and $\alpha_r(K_s)=0=\alpha_r(JK_s)$.)

The two families of $1$-forms $J\boldsymbol\theta$ and $\iota_{\boldsymbol
K}\omega$ span the same $\ell$-dimensional space. Hence we may write
\begin{equation*}
J\theta_r= G_{rs}\iota_{K_s}\omega\qquad\text{and}\qquad
\iota_{K_r}\omega= H_{rs}J\theta_s,
\end{equation*}
where $G_{rs}$ and $H_{rs}$ are mutually inverse and $\boldsymbol
K$-invariant. Since $\g{K_r,K_s}=\omega(K_r,JK_s)=H_{rs}$, we deduce that
$G_{rs}$ and $H_{rs}$ are symmetric and positive definite.

\begin{prop}\textup{\cite{ped-poon}}\label{p:ped-poon}
Let $(S,J)$ be a complex $2(m-\ell)$-manifold, let $B=S\times U$ with $U$ open
in $\R^{\ell*}$, and let $M$ be a principal $\ell$-torus bundle over
$B$. Denote the components of the projection $\boldsymbol\sigma\colon
B\to\R^{\ell*}$ by $\sigma_r$.  Now suppose that\textup:
\begin{numlist}
\item $(h,\omega_h)$ is family of compatible K\"ahler metrics on the level
surfaces of $\boldsymbol\sigma$ in $B$\textup;
\item $G_{rs}$ is a symmetric positive definite matrix of functions on
$B$, with inverse matrix $H_{rs}$, satisfying the equations
\begin{equation}\label{eq:Geqn}
\frac{\partial G_{rs}}{\partial\sigma_t}
=\frac{\partial G_{rt}}{\partial\sigma_s} \qquad\text{and}\qquad
d_hd^c_h G_{rs}
+ \frac{\partial^2\omega_h}{\partial\sigma_r\partial\sigma_s}=0;
\end{equation}
\item $\boldsymbol\theta\colon M\to \R^\ell$ is the connection $1$-form of a
principal connection on $M$ over $B$ whose curvature satisfies the equation
\begin{equation}\label{eq:thetaeqn}
d\theta_r=\frac{\partial\omega_h}{\partial\sigma_r}+d^c_hG_{rs}\wedge
d\sigma_s.
\end{equation}
\end{numlist}
Then the almost hermitian structure
\begin{equation}\label{eq:htks}\begin{split}
g&= h + G_{rs} d\sigma_r d\sigma_s + H_{rs} \theta_r\theta_s\\
\omega&= \omega_h+d\sigma_r\wedge\theta_r\\
J\theta_r& = -G_{rs}d\sigma_s
\end{split}\end{equation}
on $M$ is K\"ahler with a free hamiltonian $\ell$-torus action and
K\"ahler quotient $S$.

Any K\"ahler manifold with a local hamiltonian $\ell$-torus action arises
locally in this way on the dense open set $M^0$ where the Killing vector
fields are independent.
\end{prop}
\begin{proof} We have seen already that any K\"ahler structure with a
local hamiltonian $\ell$-torus action can be written in the
form~\eqref{eq:htks}, where $(h,J,\omega_h)$ is K\"ahler for each fixed
$\boldsymbol\sigma$ and $d\theta_r(JK_s,JK_t)=-\theta_r([JK_s,JK_t])=0$.  Now
under these conditions, $\zeta_r=G_{rs} d\sigma_s + i \theta_r$ generate the
$(1,0)$-forms on the fibres, so $(g,J,\omega)$ is K\"ahler if and only if the
$\zeta_r$ generate a differential ideal modulo horizontal forms and $\omega$
is closed. Since $d\omega=(\partial\omega_h/\partial\sigma_r -
d\theta_r)\wedge d\sigma_r$, and $d\zeta_r= d_h G_{rs} \wedge d\sigma_s -
(\partial G_{rs}/\partial\sigma_t)d\sigma_s\wedge d\sigma_t+id\theta_r$, it
follows easily that $(g,J,\omega)$ is K\"ahler if and only if $G_{rs}$
satisfies~\eqref{eq:Geqn} and $\theta_r$ satisfies~\eqref{eq:thetaeqn}: the
second part of \eqref{eq:Geqn} follows from the integrability
of~\eqref{eq:thetaeqn}.
\end{proof}

\begin{rem} It is crucial here that the local torus action
is hamiltonian in the strong sense that the components of $\boldsymbol\sigma$
Poisson-commute, i.e., $d\sigma_r(K_s)=0$. This condition is often missed in
the literature, since if $K_r$ and $K_s$ commute, $d\sigma_r(K_s)$ is
constant, and on a compact manifold $\sigma_r$ must have a critical point, so
the constant is zero. However, we are not assuming compactness: indeed, the
above local description is only valid on an open set where $\sigma_r$ has no
critical points!
\end{rem}

Proposition~\ref{p:ped-poon} shows that a K\"ahler metric with a local
hamiltonian $\ell$-torus action may be specified by essentially free
data. Indeed~\eqref{eq:thetaeqn} is integrable by virtue of~\eqref{eq:Geqn}.
To solve the latter, observe that the first part implies we can write
$G_{rs}=\partial u_r/\partial\sigma_s$, and since $G_{rs}$ is symmetric,
$u_r=\partial G/\partial\sigma_r$ for some function $G$ on $B$ such that the
$d_h$-closed $J$-invariant $2$-form $\omega_h+d_hd^c_hG$ depends affinely on
$\boldsymbol\sigma$. Since we can add a $d_hd^c_h$ potential for this $2$-form
to $G$ without altering $G_{rs}$, we can assume (locally) that for each fixed
$\boldsymbol\sigma$, $\omega_h+d_hd_h^cG=0$.  Thus $G$ determines $\omega_h$
and $G_{rs}$ and is now subject only to the open condition that these are
positive definite.

In fact $G$ is a fibrewise Legendre transform of a K\"ahler potential,
generalizing work of Guillemin in the toric case~\cite{G:kstv} (see
also~\cite{Abreu}). Observe first that $d^c u_r=d_h^c u_r+\theta_r$, and so
$dd^c u_r =\partial\omega_h/\partial\sigma_r+d_hd^c_h u_r=0$, i.e., $u_r$ is
pluriharmonic.  Now $dd^c (\sigma_ru_r-G)=d(\sigma_r\wedge d^cu_r-
d_h^cG)=d\sigma_r\wedge\theta_r-d_hd_h^cG = \omega$, so $H:=\sigma_ru_r-G$ is
a K\"ahler potential. Since $G_{rs}$ is nondegenerate, the $u_r$ also form a
coordinate system on each fibre of $B$ over $S$, and we let $\partial/\partial
u_r=H_{rs}\,\partial/\partial\sigma_s$ be the coordinate vector fields tangent
to the fibres, so that $H_{rs}=\partial\sigma_s/\partial u_r$ and
$\sigma_s=\partial H/\partial u_r$. If we locally set $d^cu_r=dt_r$ then
$\boldsymbol u+i\boldsymbol t\colon M\to\C^{\ell}$ is holomorphic and
$\theta_r=dt_r+\alpha_r$ with $\alpha_r=-d_h^c u_r$. This is the fibrewise
Legendre dual coordinate system to $(\boldsymbol\sigma,\boldsymbol t)\colon
M\to\R^{\ell*}\times\R^\ell$, and we refer to $G$ as a {\it dual potential}.

It is convenient to introduce a fixed ($\boldsymbol\sigma$-independent)
volume form $\vol_S$ on $S$ and write $\vol_{\omega_h}=Q\vol_S$. Observe in
particular that
\begin{equation}\label{eq:tromdom}
\G{\omega_h,\frac{\partial\omega_h}{\partial\sigma_r}}_{\!h}
=\G{\omega_h^{-1},\frac{\partial\omega_h}{\partial\sigma_r}}
=\frac1Q\frac{\partial Q}{\partial\sigma_r}.
\end{equation}

\begin{prop}\label{p:pluriricci} Let $(M,g,J,\omega)$ be K\"ahler with a local
hamiltonian $\ell$-torus action.
\begin{numlist}
\item Let $f$ be any invariant function on $M$. Then
\begin{align}\notag
dd^cf &= d_h d^c_h f + \frac{\partial f}{\partial u_r}
\frac{\partial\omega_h}{\partial\sigma_r}\\
\label{eq:ddcf}
&\quad+d_h\Bigl(\frac{\partial f}{\partial u_r}\Bigr)\wedge\theta_r
+d^c_h\Bigl(\frac{\partial f}{\partial u_r}\Bigr)\wedge J\theta_r\\
\notag&\quad
+\frac{\partial}{\partial\sigma_r}\Bigl( \frac{\partial f}{\partial u_s}\Bigr)
d\sigma_r\wedge\theta_s,\\
\label{eq:lapf}
\Delta f&=\Delta_h f-\frac1Q\frac{\partial}{\partial\sigma_r}
\Bigl(Q\frac{\partial f}{\partial u_r}\Bigr).
\end{align}
It follows that $dd^cf=0$ if and only if $f=A_r u_r+B$ where the $A_r$ are
constant, and $B$ is a pluriharmonic function on $S$.

\item Suppose that $\kappa_h$ is a Ricci potential for $\omega_h$ for each
fixed $\boldsymbol\sigma$, i.e., $d_hd^c_h\kappa_h=\rho_h$ where $\rho_h$ is
the Ricci form of $\omega_h$.  Then $\kappa=\kappa_h+\frac12\log\det G_{rs}$
is a Ricci potential for $g$, and we have
\begin{align}\label{eq:ddckap}
d_hd^c_h\kappa &= \rho_h+\frac12 d_h(H_{rs}d_h^c G_{rs}),\\
\frac{\partial\kappa}{\partial u_r}&=-\frac1{2Q}\frac{\partial}
{\partial\sigma_t}(QH_{rt}).\label{eq:dkapdu}
\end{align}
\end{numlist}
\end{prop}
\begin{numlproof}
\item Expanding $d$ into horizontal and vertical parts, we get
\begin{equation*}
dd^cf=d_h d_h^c f+\frac{\partial f}{\partial u_r} d\theta_r
+d\sigma_r\wedge\frac{\partial}{\partial\sigma_r}d^c_hf+
d_h\frac{\partial f}{\partial u_r}\wedge\theta_r+
\frac{\partial}{\partial\sigma_s}\frac{\partial f}{\partial u_r}\,
d\sigma_s\wedge\theta_r.
\end{equation*}
The second term in the equation~\eqref{eq:thetaeqn} for $d\theta_r$ combines
with the third term in the above equation to give $d^c_h(\partial f/\partial
u_r)\wedge J\theta_r$. The formula for the laplacian~\eqref{eq:lapf} follows
by contracting with $\omega$, using~\eqref{eq:tromdom}.

Now $dd^cf=0$ if and only if the three lines on the right hand side
of~\eqref{eq:ddcf} are separately zero.  Hence $\partial f/\partial u_r$ must
be constant, i.e., $f=A_r u_r+B$ with $A_r$ constant and $\partial B/\partial
u_r=0$, and so $dd^cf=0$ if and only if $dd^cB=d_hd_h^cB=0$.

\item A Ricci potential has the form $-\frac12\log(\vol_\omega/\vol_J)$ where
$\vol_\omega=\frac1{m!}\omega\wedge\cdots\wedge\omega$ and $\vol_J$ is a
holomorphic volume form. We first observe (see~\cite{ped-poon,christina2})
that if $dz_\mu$ is a local frame of holomorphic $(1,0)$-forms on $S$, then
there are functions $B_{r\mu}$ such that $\sum_\mu B_{r\mu}dz_\mu
+G_{rs}d\sigma_s+i\theta_r$, together with $dz_\mu$, form a local holomorphic
frame of $M$. Since $\omega=\omega_h+d\sigma_r\wedge\theta_r$, the formula for
$\kappa$ is immediate.

Equation~\eqref{eq:ddckap} follows easily using the fact that for any matrix
valued function $A$, $d\log\det A = \trace A^{-1}dA$.  For~\eqref{eq:dkapdu}
we also note that $\kappa_h=-\frac12\log(Q\vol_S/\vol_{J_h})$, where
$\vol_{J_h}$ is a ($\boldsymbol\sigma$-independent) holomorphic volume form on
$S$, and so
\begin{align*}
\frac{\partial\kappa}{\partial u_r}&=
H_{rt}\frac{\partial\kappa}{\partial \sigma_t}
=\frac12H_{rt}\Bigl(H_{pq}\frac{\partial G_{pq}}{\partial\sigma_t}
-\frac1Q\frac{\partial Q}{\partial\sigma_t}\Bigr)
=\frac12\Bigl(H_{rt}\frac{\partial G_{tp}}{\partial\sigma_q}H_{pq}
-\frac1Q\frac{\partial Q}{\partial\sigma_t}H_{rt}\Bigr)\\
&=-\frac12\Bigl(\frac{\partial H_{rt}}{\partial\sigma_t}+
\frac1Q\frac{\partial Q}{\partial\sigma_t}H_{rt}\Bigr),
\end{align*}
where we use the symmetry of $\partial G_{pq}/\partial\sigma_t$ in $p,q,t$.
\end{numlproof}

This result provides conditions for $M$ to be K\"ahler--Einstein, using the
fact that invariant Ricci and K\"ahler potentials then differ by an invariant
pluriharmonic function. More generally, substituting~\eqref{eq:ddckap}
and~\eqref{eq:dkapdu} into~\eqref{eq:ddcf} and~\eqref{eq:lapf} gives the Ricci
form and scalar curvature. These expressions are rather complicated in
general. However, if we suppose that $d_hG_{rs}=0$ and $d_hQ=0$ then
\begin{equation*}\begin{split}
\rho&= \rho_h-\frac 1{2Q}\frac{\partial(QH_{rt})}
{\partial\sigma_t}\omega_h-\frac{\partial}{\partial\sigma_s}
\Bigl(\frac1{2Q}\frac{\partial(Q H_{rt})}
{\partial\sigma_t}\Bigr)d\sigma_s\wedge\theta_r,\\
\Scal&=\Scal_h
-\frac1Q \frac{\partial^2}{\partial\sigma_r\partial\sigma_s}(QH_{rs}).
\end{split}\end{equation*}
Note that these expressions depend linearly in $H_{rs}$: this fact was
emphasised by Abreu~\cite{Abreu} in the toric case, and by
Hwang--Singer~\cite{hwang-singer} in the case of circle symmetry: we have just
combined their arguments. We shall see the significance of the conditions
$d_hG_{rs}=0$ and $d_hQ=0$ shortly.

We remark that when $m=2$, $\ell=1$, the Pedersen--Poon construction reduces
to LeBrun's construction~\cite{lebrun}:
\begin{equation*}
g=we^u(dx^2+dy^2)+w\,dz^2+w^{-1}\theta^2,\qquad
\omega=we^u\,dx\wedge dy+dz\wedge\theta,
\end{equation*}
where $w_{xx}+w_{yy}+(we^u)_{zz}=0$, which is the integrability condition for
$d\theta=w_x \,dy\wedge dz-w_y\,dx\wedge dz+(we^u)_z\,dx\wedge dy$.  Here $d_h
w= w_x \,dx+w_y dy$. Note $\kappa=-\frac12u$ is a Ricci potential for $g$.

\subsection{Rigid hamiltonian torus actions}\label{s:rigid}

K\"ahler manifolds with a hamiltonian $\ell$-torus action are too complicated,
in their fullest generality, for constructing interesting K\"ahler
metrics. Indeed most applications, including those in~\cite{ped-poon}, use
only K\"ahler metrics in the following subclass.

\begin{prop}\label{p:gcc}
Suppose the K\"ahler manifold $(M,g,J,\omega)$ has a local
\textup(isometric\textup) hamiltonian $\ell$-torus action $\boldsymbol
K=J\grad_{\smash g} \boldsymbol\sigma$, for $\boldsymbol\sigma\colon
M\to\R^{\ell*}$, and let $\cF$ be the foliation generated by $K_r, JK_r$
$(r=1,\ldots\ell)$. Then on the open dense set $M^0$, where the action is
locally free, the following are equivalent\textup:
\begin{numlist}
\item the leaves of $\cF$ are totally geodesic\textup;
\item the connection $\boldsymbol\theta\colon TM^0\to \R^\ell$, with
$\ker\boldsymbol\theta=({\rm span}\,\boldsymbol K)^\perp$, is $J\boldsymbol
K$-invariant\textup;
\item $\g{K_r,K_s}$ is constant on the level surfaces of $\boldsymbol\sigma$
for all $r,s$\textup;
\item the family of K\"ahler forms $\omega_h=\omega-d\sigma_r\wedge\theta_r$
on the local leaf space of $\cF$ depends affinely on $\boldsymbol\sigma$ and
the linear part pulls back to the curvature of $\boldsymbol\theta$\textup;
\item there is a \textup(local\textup) $\boldsymbol K$-invariant K\"ahler
potential of the form $H=H_0+H_\sigma$ where $H_0$ is constant on the leaves
of $\cF$ and $H_\sigma$ is constant on the level surfaces of
$\boldsymbol\sigma$.
\end{numlist}
\end{prop}
\begin{proof}  The conditions (i)--(iii) are all equivalent to the fact that
$\g{\nabla_{K_r}K_s,X}=0$ for all $r,s$ and all $X$ orthogonal to
$\cF$. Indeed, since $J$ is parallel and $K_r$ is holomorphic this says that
$\cF$ is totally geodesic.  On the other hand, $\g{\cL_{JK_r}(JX),K_s}=
\g{J\cL_{JK_r}(X),K_s}=\g{J\nabla_{JK_r}X+\nabla_XK_r,K_s}=
-2\g{\nabla_{K_r}K_s,X}$, so it also says that the connection
$\boldsymbol\theta$ is $JK_r$-invariant. Finally, it means that
$\partial_X\g{K_r,K_s}=-2\g{\nabla_{K_r}K_s,X}=0$ for all $X$ orthogonal
to $\cF$.

To establish the equivalence of the local conditions (iii)--(v), we use the
Pedersen--Poon construction, Proposition~\ref{p:ped-poon}. (iii) means that
$d_hH_{rs}=0$, or equivalently $d_hG_{rs}=0$, which by~\eqref{eq:thetaeqn} is
equivalent to $d\theta_r=\partial\omega_h/\partial\sigma_r$; this is (iv),
since~\eqref{eq:Geqn} then shows that $\omega_h$ is affine in
$\boldsymbol\sigma$. (v) gives that $H_{rs}$ is the hessian of $H_\sigma$
which implies (iii). Conversely (iii) implies that the dual potential $G$ is
an affine function of $\boldsymbol\sigma$, so that $d_hH=
\sigma_r(\partial/\partial\sigma_r)(d_h G) - d_h G$ is independent of
$\boldsymbol\sigma$; then it has a local $\boldsymbol\sigma$-independent
$d_h$-potential $H_0$, and $d_h H_\sigma=0$ where $H_\sigma=H-H_0$.
\end{proof}

If $M$ is given by the Pedersen--Poon construction (as it is locally), then
(i) means that the fibres of $M\to S$ (the complex orbits) are totally
geodesic, (ii) that $M\to B$ is the pullback of a principal bundle with
connection over $S$, and (iii) that the metric on the fibres of $M\to B$ (the
torus orbits) depends only on the momentum map; the condition (iv) on the
K\"ahler quotient is a kind of rigid Duistermaat--Heckman property (it holds
in cohomology by~\cite{DH}), while (v) generalizes Calabi's
Ansatz~\cite{calabi0} for K\"ahler metrics on holomorphic bundles.

\begin{defn} A local hamiltonian $\ell$-torus action $\boldsymbol K=
J\grad_{\smash g} \boldsymbol\sigma$ on a K\"ahler manifold will be called
{\it rigid} if $\g{K_r,K_s}$ is constant on the level surfaces of
$\boldsymbol\sigma$.
\end{defn}

\begin{prop}\label{p:lcc} Suppose that $M$ arises from the Pedersen--Poon
construction for a rigid hamiltonian $\ell$-torus action, and let $\nabla^\|$
and $\nabla^H$ be respectively the Levi-Civita connection on the fibres of $M$
over $S$, and the Levi-Civita connection on the level surfaces of
$\boldsymbol\sigma$ in $B$, lifted to the horizontal distribution of $M\to S$.
Let $X,Y$ be horizontal vector fields and $U,V$ be vertical vector fields.
Then
\begin{align}
\label{eq:hor}&&\nabla_X Y &\quad= && \nabla^H_XY  &&- C(X,Y) &&\\
\label{eq:hv} &&\nabla_X U &\quad=
&& \hphantom{[U,X]^H+{}}\g{C(X,\cdot),U} &&+ [X,U]^{\smash\|} &&\\
\label{eq:vh}&&\nabla_U X &\quad= && [U,X]^H+\g{C(X,\cdot),U}&& &&\\
\label{eq:vert}&&\nabla_U V &\quad= && &&\nabla^{\smash[t]\|}_U V, &&
\end{align}
where ${}^H$ and ${}^\|$ denote the horizontal and vertical components,
and the O'Neill tensor $C$ is given by
\begin{equation}\label{eq:oneill}
2C(X,Y) = \Omega_r(X,Y) K_r + \Omega_r(JX,Y) JK_r.
\end{equation}
\end{prop}
\begin{proof} These observations all follow from the
Koszul formula
\begin{equation}\label{eq:koszul}\begin{split}
2\g{\nabla_XY,Z}&=\partial_X\,\g{Y,Z}+\partial_Y\,\g{X,Z}
-\partial_Z\,\g{X,Y}\\
&+\g{[X,Y],Z}-\g{[X,Z],Y}-\g{[Y,Z],X}.
\end{split}\end{equation}
The contraction of~\eqref{eq:hor} with a horizontal vector field $Z$ is
immediate because $M\to B$ is a riemannian submersion.  For the vertical
component, \eqref{eq:koszul} gives
\begin{equation*}
2\g{\nabla_XY,Z}=(\cL_Z g)(X,Y)+\g{[X,Y],Z}.
\end{equation*}
(with $Z$ vertical). Taking $Z=\partial/\partial\sigma_r$ and $Z=K_r$ we
obtain~\eqref{eq:hor} with $C$ given by~\eqref{eq:oneill}, since
$[X,Y]^\|=-\Omega_r(X,Y)K_r$, $JK_r=-\grad_g\sigma_r$ and
\begin{equation*}
\frac{\partial g}{\partial\sigma_r}(X,Y)
=-\frac{\partial\omega_h}{\partial\sigma_r}(JX,Y)=
-\Omega_r(JX,Y).
\end{equation*}
The remaining three equations are much easier: $\g{\nabla_XU,Y}= -\g{\nabla_X
Y,U}$, $\nabla_U X-\nabla_XU=[U,X]$ and $\g{\nabla_U X,V} =-\g{\nabla_U V,X}$,
so we only need to check $\g{\nabla_X U,V}=\g{[X,U],V}$
and~\eqref{eq:vert}. These follow immediately because the metric on the fibres
is constant along horizontal curves and the fibres are totally geodesic.
\end{proof}

\subsection{Semisimple K\"ahler quotients}\label{s:ss}

\begin{defn} A complex manifold $(S,J)$ with a family of K\"ahler metrics
$(h,\omega_h)$ (with parameter $\boldsymbol\sigma$) is {\it semisimple} if
there is a K\"ahler form $\Omega_S$ on $S$ with respect to which the
$\omega_h$ are simultaneously diagonalizable and parallel.
A local hamiltonian torus action is semisimple if its local K\"ahler quotient
is.
\end{defn}
We can of course take $\Omega_S$ to be $\omega_h$ for some fixed
$\boldsymbol\sigma$, but it will be convenient later to make a different
choice.
\begin{prop} If $(S,J)$ is semisimple then $(S,\Omega_S)$ is a locally
a K\"ahler product of $(S_a,\omega_a)$ $(a=1,\ldots N,\, N\geq1)$ such that
$\omega_h=\sum_{a=1}^N c_a(\boldsymbol\sigma)\omega_a$, where
$c_a(\boldsymbol\sigma)$ is constant on $S$. The Levi-Civita connection of
$\omega_h$ is independent of $\boldsymbol\sigma$, being equal to the
Levi-Civita connection of $\Omega_S$.
\end{prop}
\begin{proof}
As $\omega_h$ is a parallel $(1,1)$-form on $S$, the ($J$-invariant,
simultaneous) eigendistributions of $\omega_h$ are parallel, and $S$
splits as a local K\"ahler product by the deRham theorem. The Levi-Civita
connections of
$\omega_h$ and $\Omega_S$ agree, because on each factor $S_a$ of the
local K\"ahler product they are related by a constant multiple.
\end{proof}
In the case of local hamiltonian $\ell$-torus actions, the semisimplicity
condition implies in particular that the quantity $Q$ defined by
$\vol_{\smash{\omega_h}}=Q\vol_S$, where $\vol_S$ is the volume form of
$\Omega_S$, is constant on $S$, being given by
\begin{equation}\label{eq:Q}
Q=\textstyle\prod_{a=1}^N c_a(\boldsymbol\sigma)^{m_a},\qquad
\text{where}\qquad \dim S_a=2m_a.
\end{equation}
If the action is also rigid, $\omega_h$ depends affinely on
$\boldsymbol\sigma\in\R^{\ell*}$, so we can write
$\omega_h=\Omega_0+\langle\boldsymbol\sigma, \boldsymbol\Omega\rangle$ where
$\Omega_0$ and $\boldsymbol\Omega$ are closed $J$-invariant
($\boldsymbol\sigma$-independent) $2$-forms on $S$, the latter with values in
$\R^\ell$. Letting $\Omega_r$ denote the components of $\boldsymbol\Omega$, we
have $\Omega_r=\sum_{a=1}^N c_{ar}\omega_a$ for $r=0,\ldots \ell$, where
$c_a(\boldsymbol\sigma)=c_{a0}+c_{ar}\sigma_r$.

\subsection{Orthotoric K{\"a}hler metrics}\label{s:ortho}

\begin{defn} A K{\"a}hler $2m$-manifold $(M, g, J, \omega)$ is
{\it orthotoric} if it is equipped with $m$ Poisson-commuting functions
$\sigma_1,\ldots\sigma_m$ such that $K_r = J {\rm grad} _{\smash g} \sigma _r$
are Killing vector fields and, on a dense open set $M^0$, the roots $\xi_j$ of
$\sum_{r=0}^m (-1)^r\sigma_r t^{\ell-r}$ $(\sigma_0=1$) are smoothly defined,
with linearly independent, orthogonal gradients.
\end{defn}

Note that an orthotoric K\"ahler manifold is toric, and any toric Riemann
surface is orthotoric. For a higher dimensional toric manifold it is hard to
detect whether it is orthotoric, since the condition depends on a choice of
basis for Lie algebra of the torus. Because of this choice, we abandon the
summation convention.

The exterior derivative of the identity
\begin{equation*}
\prod _{k = 1} ^m (t-\xi _k) = \sum_{r=0}^m (-1)^r \sigma_r t^{m-r}
\end{equation*}
at $t=\xi_j$ yields
\begin{equation}\label{eq:sigtoxi}
d\xi_j=\frac1{\Delta_j}\sum_{r=0}^m (-1)^{r-1}\xi_j^{m-r}d\sigma_r,
\end{equation}
where $\Delta_j=\prod_{k\neq j}(\xi _j - \xi _k)$. This is inverse to the
identity
\begin{equation}
d\sigma_r=\sum_{j=1}^m \sigma_{r-1}(\hat\xi_j)d\xi_j,
\end{equation}
where $\sigma _{r - 1} (\hat{\xi} _j)$ denote the elementary symmetric
functions of the $m - 1$ functions $\xi _k$ with $\xi _j$ deleted (with the
convention that $\sigma _0 = 1$). Hence the coordinate systems given by
$\xi_j$ and $\sigma_r$ are related by the Vandermonde matrix and its
inverse. We have collected some Vandermonde identities that we need in
Appendix B.

\begin{prop}\label{p:ortho}
Let $(M, g, J, \omega)$ be an orthotoric K{\"a}hler $2m$-manifold.  Then, on
any simply connected domain $U$ in $M^0$, there are $m$ functions $t _r$, each
determined up to an additive constant, and $m$ functions $\oF_j$ of one
variable, such that $\{ \xi _1, \ldots \xi _m, t _1, \ldots t _m \}$ form a
coordinate system with respect to which the K{\"a}hler structure may be
written
\begin{equation} \label{eq:otks}\begin{split}
g &= \sum _{j = 1} ^m \frac{\Delta_j}
{\oF_j (\xi _j)} \, d \xi _j^2
+ \sum _{j = 1} ^m \frac{\oF_j (\xi _j)}{\Delta_j}
\Bigl (\sum _{r = 1} ^m \sigma _{r - 1} (\hat{\xi} _j) \, dt _r \Bigr )^2,\\
\omega &= \sum _{j=1} ^m d \xi _j \wedge\Bigl (
\sum _{r = 1} ^m \sigma _{r - 1} (\hat{\xi} _j) d t _r \Bigr )
= \sum _{r = 1} ^m d\sigma_r \wedge d t_r,\\
J d \xi _j &= \frac{\oF_j (\xi _j)}{\Delta_j} \,
\sum _{r = 1} ^m \sigma _{r - 1} (\hat{\xi} _j) \, d t _r, \qquad\qquad
J d t _r = (-1) ^r \,\sum_{j=1}^m \frac{\xi_j ^{m-r}}{\oF_j(\xi _j)}\,d \xi _j.
\end{split}\end{equation}

Conversely, for any $m$ real functions $\oF_j$ of one variable, the
almost-hermitian structure defined
by~\eqref{eq:otks} is K{\"a}hler and orthotoric
with dual potential
\begin{equation}\label{eq:otG}
G=-\sum_{j=1}^m\int^{\xi_j} \frac{\prod_k(t-\xi_k)}{\oF_j(t)}dt
=-\sum_{j=1}^m\int^{\xi_j} \frac{\sum_{r=0}^m(-1)^r\sigma_r t^{m-r}}
{\oF_j(t)}dt
\end{equation}
and K\"ahler potential
\begin{equation}\label{eq:otH}
H=\sum_{j=1}^m\int^{\xi_j} \frac{t^m}{\oF_j(t)}dt.
\end{equation}
\end{prop}
\begin{proof} We apply Proposition~\ref{p:ped-poon} to obtain the local
expression.
The condition that the $\xi_j$ have orthogonal gradients means that
\begin{equation}
H_{rs}=\sum_{j=1}^m \sigma_{r-1}(\hat\xi_j)\sigma_{s-1}(\hat\xi_j)
|d\xi_j|^2
\end{equation}
and hence
\begin{equation}
G_{rs}=\sum_{j=1}^m \frac{(-1)^{r+s}\xi_j^{m-r}\xi_j^{m-s}}{\Delta_j^2
|d\xi_j|^2}.
\end{equation}
We set $\oF_j=\Delta_j|d\xi_j|^2$ so that
\begin{equation}\label{eq:dur}
\sum_{s=1}^m G_{rs}d\sigma_s=\sum_{j=1}^m
\frac{(-1)^{r-1}\xi_j^{m-r}}{\oF_j}d\xi_j.
\end{equation}
If $G_{rs}$ is the hessian of a function $G$, this must be closed,
i.e.,
\begin{equation*}
\xi _j ^{m - r} \frac{\partial \oF_j}{\partial \xi _k}
= \xi _k ^{m - r} \frac{\partial \oF_k}{\partial \xi _j},
\end{equation*}
for all $j, k, r = 1, \ldots m$. Multiplying this by $(-1)^rd\sigma_r$ and
using \eqref{eq:sigtoxi} to sum over $r$ (which amounts to inverting the above
Vandermonde system) shows that $\Delta_j\partial \oF_j/\partial\xi_k$ vanishes
for $j\neq k$; it follows that each $\oF_j$ only depends on $\xi _j$.
Changing coordinates from Proposition~\ref{p:ped-poon} yields~\eqref{eq:otks}.

Conversely if $\oF_j$ is a function only of $\xi_j$ then~\eqref{eq:dur}
is equal to $du_r$ where
\begin{equation}
u_r=-\sum_{j=1}^m\int^{\xi_j} \frac{(-1)^{r}t^{m-r}}{\oF_j(t)}dt.
\end{equation}
Since the integrand in~\eqref{eq:otG} vanishes when $t=\xi_j$, the derivative
of $G$ with respect to $\sigma_r$ is $u_r$, so its hessian is $G_{rs}$;
$\sum_{r=1}^m u_r\sigma_r-G$ then gives~\eqref{eq:otH}.  Since the $d \xi_j$
are evidently pairwise orthogonal, the structure is orthotoric.
\end{proof}

We end with an alternative characterization of orthotoric K\"ahler metrics.

\begin{prop}\label{p:otchar} A toric K\"ahler structure $(g,J,\omega)$ is
orthotoric if and only if there is a momentum map $(\sigma_1,\ldots \sigma_m)$
such that $\sum_{r=1}^m (\sigma_1\sigma_r-\sigma_{r+1})du_r$ is a closed
$1$-form, where $u_r=\partial G/\partial\sigma_r$ for a dual potential $G$.
\end{prop}
\begin{proof}
Let $\xi_j$ be the roots of the polynomial $\sum_{r=0}^m (-1)^r\sigma_r
t^{m-r}$.  Then
\begin{align*}
\sum_{r=1}^m d(\sigma_1\sigma_r-\sigma_{r+1})\wedge du_r
&=\sum_{r,s} G_{rs} (\sigma_r\,d\sigma_1-d\sigma_{r+1})\wedge d\sigma_s\\
&=\sum_{r,s,j,k} G_{rs}(\xi_j\,\sigma_{r-1}(\hat\xi_j)d\xi_j)\wedge
(\sigma_{s-1}(\hat\xi_k)d\xi_k),
\end{align*}
which is zero if and only if $\sum_{r,s}
G_{rs}(\xi_j-\xi_k)\sigma_{r-1}(\hat\xi_j)\sigma_{s-1}(\hat\xi_k)=0$
for all $j,k$. The left hand side is
$(\xi_j-\xi_k)\g{\partial/\partial\xi_j, \partial/\partial\xi_k}$, so
the result follows.
\end{proof}

\section{Classification of hamiltonian $2$-forms}\label{s:chf}

\subsection{Rough classification of hamiltonian $2$-forms}\label{s:rough}

On any K\"ahler manifold, any $J$-invariant parallel $2$-form is hamiltonian.
However, in this case the Killing vector fields $K(t)$ are all identically
zero. For a general hamiltonian $2$-form $\phi$, it is important to know how
many of the $K(t)$ are linearly independent. To do this, we (temporarily)
write $\Mp(t)=\sum_{r=0}^m (-1)^r \sigma_r t^{m-r}$ and $K_r=J\grad_g
\sigma_r$ for the coefficients of the momentum polynomial and corresponding
Killing vector fields. Hence $K(t)=\sum_{r=1}^m (-1)^r K_r t^{m-r}$ is a
linear combination of $K_1,\ldots K_m$ for any $t$.

\begin{prop} Let $\phi$ be a hamiltonian $2$-form on a
\textup(connected\textup) K\"ahler $2m$-manifold $M$. Then there is an integer
$\ell$, with $0\leq \ell\leq m$ such that $K_1\wedge\cdots\wedge K_\ell$ is
nonzero on a dense open subset, but $\dim{\rm span}\{K_1,\ldots K_m\}\leq\ell$
on all of $M$.
\end{prop}
\begin{proof} The coefficient of $t^{m-r}$ in identity~\eqref{Krelation}
gives
\begin{equation}\label{Krelr}
K_{r+1}=\phi(JK_r)+\sigma_r K_1.
\end{equation}
Suppose for some $z\in M$ and $1< s < m$, $K_s$ is a linear combination of
$K_1,\ldots K_{s-1}$ at $z$.  Then $K_{s+1}$ is also a linear combination of
$K_1,\ldots K_{s-1}$ at $z$: to see this, use~\eqref{Krelr} with $r=s$, then
write $\phi(JK_s)$ as a linear combination of $\phi(JK_r)$ with $r<s$ and
use~\eqref{Krelr} again to express these in terms of $K_1,\ldots K_s$. Hence,
at each $z\in M$, $\dim{\rm span}\{K_1,\ldots K_m\}$ is the largest integer
$\ell_z$ such that $K_1,\ldots K_{\ell_z}$ are linearly independent at
$z$. However, for any integer $r$, $K_1,\ldots K_r$ are linearly dependent if
and only if the holomorphic $r$-vector $K_1^{1,0}\wedge\cdots\wedge K_r^{1,0}$
is zero.  Hence the set where $K_1,\ldots K_r$ are linearly independent is
empty or dense. The result follows.
\end{proof}

The integer $\ell$ of this proposition will be called the \emph{order} of
$\phi$ and we let $M^0$ be the dense open set where $K_1,\ldots K_\ell$ are
independent. We shall identify the order of $\phi$ with the number of
non-constant roots of the momentum polynomial $\Mp$.

\begin{lemma}\label{lem:evec} If, on an open subset of $M$,
$\phi (Z,\cdot) = \xi\, \omega(Z,\cdot)$ with $Z$ nonvanishing, then $d\xi$ is
the orthogonal projection of $d\sigma$ onto the span of $Z$ and $JZ$.
\end{lemma}
\begin{proof} Without loss of generality, we can take $Z$ to be a unit
vector field, and hence $(\nabla _X \phi) (Z, JZ) = d \xi (X)$ for all vector
fields $X$. By \eqref{ham}, this becomes
\begin{equation}\label{projdxi}
d \xi = d \sigma (Z) Z + d \sigma (JZ) JZ,
\end{equation}
which is what we wanted to prove.
\end{proof}

The roots of the momentum polynomial are the eigenvalues of
$-J\phi=-J\circ\phi$, viewed as a $J$-commuting symmetric endomorphism of
$TM$. At each point of $M$, these eigenvalues are real and there is an
orthogonal $J$-invariant direct sum decomposition of the tangent space into
eigenspaces. We count an eigenvalue with multiplicity $k$ if the corresponding
eigenspace has real dimension $2k$; for the moment, we denote by $\xi _1,
\ldots \xi _m$ the $m$ (not necessarily distinct) eigenvalues of $-J\phi$. It
follows that, for any $t$,
\begin{equation*}
\Mp(t) = \prod _{j = 1} ^m (t-\xi _j)
= t^m - \sigma _1 t^{m - 1} + \cdots + (- 1) ^m \sigma _m,
\end{equation*}
where $\sigma_1,\ldots \sigma_m$ are the elementary symmetric functions of
$\xi _1,\ldots \xi_m$.  The above lemma, and the independence of $K_1,\ldots
K_\ell$ on $M^0$, yields a fundamental fact.

\begin{prop}\label{p:roots} Let $\phi$ be a hamiltonian $2$-form on $M$.
Then the roots $\xi_i$ of $\Mp(t)$ and their derivatives $d\xi_i$ may be
defined smoothly on the dense open set $M^0$ and the roots extend
continuously to $M$.  Furthermore, for $i\neq j$, $d\xi_i$ and $d\xi_j$ are
orthogonal on $M$.  In particular, any repeated root \textup(on an open
set\textup) is constant.

If the order of $\phi$ is $\ell$, then there are $\ell$ non-constant roots
and they are functionally independent on the dense open set where
$K_1,\ldots K_\ell$ are independent.
\end{prop}
\begin{proof} The ordered roots $\xi_1\leq \cdots\leq \xi_m$ are
continuous on $M$. By the maximality of $\ell$, wherever $\xi_1,\ldots \xi_m$
can be smoothly defined, \emph{at most} $\ell$ of the $d\xi_i$ are
independent. On the other hand, on any open subset of $M^0$ where
$\xi_1,\ldots \xi_m$ are smoothly defined, \emph{at least} $\ell$ of the
$d\xi_i$ are independent. It follows that the $\xi_i$ and $d\xi_i$ can be
defined smoothly on $M^0$, that precisely $\ell$ of the $\xi_i$ are
functionally independent there. Lemma~\ref{lem:evec} now shows that for $i\neq
j$, $d\xi_i$ and $d\xi_j$ are orthogonal.
\end{proof}

The order expresses the extent to which $\phi$ constrains the K\"ahler
geometry of $M$. At one extreme, when $\ell=0$, we have an orthogonal
$J$-invariant decomposition of $TM$ into eigenspaces of $-J\phi$, and it is
easy to see that this makes $M$ into a K\"ahler product: if we write
$\omega=\sum_\xi \omega_\xi$ and $\phi=\sum_\xi \xi\omega_\xi$ where
$\omega_\xi$ is the restriction of the K\"ahler form to the $\xi$ eigenspace,
then the closedness of $\phi$ and $\omega$ is equivalent to $d\omega_\xi=0$
for each $\xi$---for instance there could be only one eigenspace (in which
case $\phi$ is a constant multiple of the K\"ahler form), or there could be
$m$ (in which case $M$ is an arbitrary K\"ahler product of Riemann surfaces).
At the other extreme, when $\ell=m$, $M$ is toric and Lemma~\ref{lem:evec}
shows in fact that it is orthotoric, so has the explicit form of
Proposition~\ref{p:ortho}, determined by $m$ functions of one variable.

\subsection{Explicit description of the metric}

We now present the general description of K\"ahler metrics with a
hamiltonian $2$-form. To do this, it will be convenient to adopt different
notation from the previous subsection.

\begin{defn} Let $\Mp(t)=(-1)^m\pfaff(\phi-t\omega)$ be the momentum
polynomial of a hamiltonian $2$-form $\phi$ of order $\ell$, and let
$\xi_1,\ldots\xi_\ell$ the \emph{non-constant} roots of $\Mp$. We denote by
$\sigma_0,\ldots\sigma_\ell$ the elementary symmetric functions of these
non-constant roots, set $K_r:=J\grad_g \sigma_r$ for $r=1,\ldots\ell$, and
write $\Mp(t)=\Mpc(t)\Mpn(t)$, where
\begin{equation}\begin{split}
\Mpn(t)&:=\prod _{j = 1} ^\ell (t-\xi _j)
= t^\ell - \sigma_1 t^{\ell-1} + \cdots + (-1)^\ell \sigma_\ell,\\
\Mpc(t)&:=\prod\nolimits_\xi (t-\xi)^{m_\xi},
\end{split}\end{equation}
and the product over $\xi$ denotes the product over the different
\emph{constant} roots of $\Mp(t)$, $m_\xi$ being the multiplicity of
the root $\xi$.
\end{defn}
If we denote the $\sigma$'s of the previous subsection by
$\tilde\sigma_1,\ldots\tilde\sigma_m$, then for $\ell=m$
$\tilde\sigma_r=\sigma_r$; otherwise, since the roots of $\Mpc$ are
constant, it follows that $d\tilde\sigma_1,\ldots d\tilde\sigma_m$ are
constant linear combinations of $d\sigma_1,\ldots d\sigma_\ell$. Now
$d\tilde\sigma_1,\ldots d\tilde\sigma_\ell$ are linearly independent. Hence
$d\sigma_1,\ldots d\sigma_\ell$ are also linearly independent, and
$\sigma_1,\ldots\sigma_\ell$ are constant affine linear combinations of
$\tilde\sigma_1,\ldots\tilde\sigma_\ell$. (Note in particular that
$d\sigma_1=d\tilde\sigma_1=d\sigma$.)

Hence $K_1,\ldots K_\ell$ generate a local hamiltonian $\ell$-torus action
and $M^0$ is locally a bundle over the K\"ahler quotient $S$ with toric
fibres.  The tangent space ${\mathcal V}$ to the fibres of $M$ over $S$ is
spanned by $K_1,\ldots K_\ell$ and $J K_1,\ldots J K_\ell$, while the
orthogonal distribution $\mathcal H$ is the direct sum of the
eigendistributions ${\mathcal H}_\xi$ corresponding to the constant
eigenvalues $\xi$ of $-J\phi$. We let $\Omega$ be the the ${\mathcal
V}$-valued $2$-form on ${\mathcal H}$ defined by
$\Omega(Y,Z)=[Y,Z]^{\mathcal V}$, the orthogonal projection of the Lie
bracket onto $\mathcal V$.

\begin{lemma}\label{l:chf}\begin{numlist}
\item For all $\xi$ and $r$, the distribution ${\mathcal H}_\xi$ is $K_r$
and $JK_r$ invariant, and descends to a parallel distribution on $S$
\textup(with respect to each quotient metric\textup).
\item If $Y$ and $Z$ belong to ${\mathcal H}_\xi$ and ${\mathcal H}_\eta$
for distinct constant eigenvalues $\xi\neq\eta$ then $\Omega(Y,Z)=0$.  If
instead $Y$ and $Z$ both belong to ${\mathcal H}_\xi$, then
\begin{equation}\label{frob}
\Omega(Y,Z) = \frac{\omega(Y,Z)}{\Mpn(\xi)}
\sum_{r=1}^\ell (-1)^{r-1}\xi^{\ell-r} K_r.
\end{equation}
\end{numlist}
\end{lemma}
\begin{proof}
Suppose that $\phi(Z)=\xi JZ$ for a constant root $\xi$. Then
\begin{equation*}
\phi(\nabla_Y Z)=\nabla_Y(\xi JZ)-(\nabla_Y\phi)(Z)
=\xi J\nabla_Y Z - \tfrac12\iota_Z (d\sigma_1\wedge JY-d^c\sigma_1\wedge Y)
\end{equation*}
for any vector field $Y$, and hence
\begin{equation}\label{phinabla}
2(\phi-\xi J)(\nabla_Y Z)=\omega(Y,Z)d\sigma_1-\g{Y,Z}Jd\sigma_1
\end{equation}
since $d\sigma_1(Z)=0=Jd\sigma_1(Z)$ (cf.~Lemma~\ref{lem:evec} and
equation~\eqref{projdxi}).

\begin{numlist}
\item We apply~\eqref{phinabla} with $Z$ in ${\mathcal H}_\xi$ and $Y$
orthogonal to ${\mathcal H}_\xi$ to deduce that $\nabla_Y Z$ also belongs to
${\mathcal H}_\xi$. Now $\cL_Y Z =\nabla_Y Z - \nabla_Z Y$, and the first
term is in ${\mathcal H}_\xi$ for $Y=K_r$ or $JK_r$. On the other hand, for
$X$ orthogonal to ${\mathcal H}_\xi$, $\g{\nabla_Z JK_r,JX} =
\g{\nabla_ZK_r,X} = -\g{\nabla_X K_r,Z}=\g{K_r,\nabla_X Z}=0$, so that
$\nabla_Z Y$ is also in ${\mathcal H}_\xi$ for $Y=K_r$ or $JK_r$.  Thus
${\mathcal H}_\xi$ descends to $S$, and since the Levi-Civita connections of
the K\"ahler quotient metrics lift to the horizontal $(\mathcal H)$ part of
$\nabla$, this distribution is parallel by~\eqref{phinabla}.

\item This again follows from~\eqref{phinabla}: if $Y$ and $Z$ belong to
distinct eigenspaces, then $\omega(Y,Z)=\g{Y,Z}=0$, so that $[Y,Z]$ is in
$\mathcal H$; otherwise, if they both belong the the $\xi$ eigenspace, we
have
\begin{align*}\notag
(\phi-\xi J)([Y,Z])&=\omega(Y,Z)d\sigma_1=\omega(Y,Z)\sum_{j=1}^\ell d\xi_j
=\omega(Y,Z)\sum_{j=1}^\ell\frac{(\phi-\xi J) Jd\xi_j}{\xi-\xi_j}.
\displaybreak[0]\\
\tag*{Hence}
\Omega(Y,Z)&=\omega(Y,Z)\sum_{j=1}^\ell\frac{Jd\xi_j^\sharp}{\xi-\xi_j}
=\frac{\omega(Y,Z)}{\prod_{k=1}^\ell (\xi-\xi_k)}
\sum_{j=1}^\ell\biggl(\prod_{k\neq j}(\xi-\xi_k) \biggr)Jd\xi_j^\sharp
\end{align*}
from which~\eqref{frob} easily follows, since $K_r=Jd\sigma_r^\sharp
=\sum_{j=1}^\ell\sigma_{r-1}(\hat\xi_j)Jd\xi_j^\sharp$.
\endnumlproof

This lemma, with Propositions~\ref{p:gcc}, \ref{p:ortho} and~\ref{p:roots},
yields our classification.

\begin{thm}\label{t:one}
Let $(M,g,J,\omega)$ be a connected K\"ahler $2m$-manifold with a hamiltonian
$2$-form $\phi$ of order $\ell$. Then there are functions
$\pF_1,\ldots\pF_\ell$ of one variable such that on a dense open subset $M^0$,
the K\"ahler structure may be written
\begin{equation}\label{eq:explicit}\begin{split}
g&=\sum_\xi \Mpn(\xi)g_{\xi}
+\sum_{j=1}^\ell \frac{\Mp'(\xi_j)}{\pF_j(\xi_j)} d\xi_j^2
+\sum_{j=1}^\ell \frac{\pF_j(\xi_j)}{\Mp'(\xi_j)}\Bigl(\sum_{r=1}^\ell
\sigma_{r-1}(\hat\xi_j)\theta_r\Bigr)^2,\\
\omega&=\sum_\xi \Mpn(\xi)\omega_{\xi}
+\sum_{r=1}^\ell d\sigma_r\wedge \theta_r,\qquad\qquad\,
d\theta_r=\sum_\xi(-1)^r\xi^{\ell-r}\omega_{\xi},\\
J d \xi _j &= \frac{\pF_j (\xi _j)}{\Mp'(\xi_j)} \,
\sum _{r = 1} ^\ell \sigma _{r - 1} (\hat{\xi} _j) \,\theta_r,\qquad\qquad
\qquad J \theta_r = (-1) ^r
\,\sum_{j=1}^\ell \frac{\Mpc(\xi_j)}{\pF_j(\xi _j)}\xi_j^{\ell-r}\,d \xi _j.
\end{split}\end{equation}
where summation over $\xi$ denotes the sum over the different constant roots
of $\Mp(t)$, $\sigma_{r-1} (\smash{\hat{\xi}_j})$ denote the elementary
symmetric functions of the $\ell - 1$ functions $\xi_k$ with $\xi_j$ deleted,
$\Mp'(t)$ is the derivative of the momentum polynomial $\Mp(t)$ with respect
to $t$, and $\pm(g_\xi,\omega_\xi)$ is a K\"ahler metric on a manifold $S_\xi$
of the same dimension as the $\xi$-eigenspace of $-J\phi$.  Dual and K\"ahler
potentials for $(g,J,\omega)$ are given by
\begin{align}\label{eq:G}
G&=-\sum_{r=0}^\ell H_r\sigma_r-\sum_{j=1}^\ell\int^{\xi_j}
\frac{\Mp(t)}{\pF_j(t)}dt,\\
H&=H_0+\sum_{j=1}^\ell\int^{\xi_j} \frac{\Mpc(t) t^\ell}{\pF_j(t)}dt,
\label{eq:H}
\end{align}
where $H_r$ is a \textup($\boldsymbol\sigma$-independent\textup) $dd^c$
potential for $\Omega_r=\sum_\xi(-1)^r\xi^{\ell-r}\omega_{\xi}$.

Furthermore, in these coordinates, the hamiltonian $2$-form may be
written
\begin{equation}\label{eq:ham?}\begin{split}
\phi &= \sum_{\smash\xi} \xi\, \Mpn(\xi)
\omega_\xi
+ \sum_{j=1}^\ell \xi_j
d\xi_j\wedge\Bigl(\sum_{r=1}^\ell\sigma_{r-1}(\hat\xi_j)\theta_r\Bigr)\\
&= \sum_{\smash\xi} \sum_{r=0}^\ell (-1)^r\sigma_r \xi^{\ell+1-r}\omega_\xi
+\sum_{r=1}^\ell (\sigma_r d\sigma_1 -d\sigma_{r+1})\wedge\theta_r.
\end{split}\end{equation}
We also have a local $dd^c$ potential for the closed form
$\phi+\sigma_1\omega$\textup:
\begin{equation}
\phi+\sigma_1\omega=\sum_\xi \xi^{\ell+1}\omega_\xi+
\sum_{r=1}^\ell d\bigl((\sigma_1\sigma_r-\sigma_{r+1})\theta_r\bigr)
=dd^c\Phi
\end{equation}
where
\begin{equation}
\Phi=-H_{-1}+\sum_{j=1}^\ell\int^{\xi_j}\frac{\Mpc(t) t^{\ell+1}}{\pF_j(t)}dt
\end{equation}
and $H_{-1}$ is a \textup($\boldsymbol\sigma$-independent\textup) $dd^c$
potential for $\Omega_{-1}=-\sum_\xi \xi^{\ell+1}\omega_\xi$.
\end{thm}
\begin{proof} By Lemma~\ref{l:chf}, the distribution $\mathcal H$ is
preserved by $JK_1,\ldots JK_\ell$; hence, by Proposition~\ref{p:gcc}, the
local fibration of $M^0$ over $S$ is totally geodesic, and the toric
structure on the fibres is constant on the level surfaces of
$\boldsymbol\sigma$.  Therefore, the restriction of $\phi$ to any fibre is
hamiltonian (even when $\ell=1$ since the trace of $\phi$ is a hamiltonian
for a Killing vector field tangent to the fibres). By
Proposition~\ref{p:roots}, the fibres are orthotoric, and the functions
$G_{rs}$ are independent of the fibre, hence so are the functions
$\oF_j=\Delta_j|d\xi_j|^2$ defining the orthotoric structure in
Proposition~\ref{p:ortho}.

Again using Lemma~\ref{l:chf}, for the constant eigenvalues $\xi$,
${\mathcal H}_{\xi}$ descends to a $J$-invariant distribution on $S$, and
$TS$ is the direct sum of these distributions. For each fixed value of
$\boldsymbol\sigma$, these distributions are parallel with respect to the
K\"ahler quotient metric $(h,\omega_h)$ and so $S$ splits locally as a
K\"ahler product of manifolds $S_\xi$.  Furthermore, the curvature $\Omega$
of $\mathcal H$ descends to $S$, so that the $2$-form
$\omega_\xi=\omega/\Mpn(\xi)$, appearing in the formula~\eqref{frob},
descends to give a K\"ahler structure on $S_\xi$, after restricting it to
the $\xi$-eigenspace distribution. (Note however, that this K\"ahler
structure will be negative definite if $\Mpn(\xi)$ is negative.)

Clearly $\omega_h=\sum_\xi \Mpn(\xi)\omega_\xi
=\sum_{r=0}^\ell\sigma_r\Omega_r$ with
\begin{equation}\label{eq:Omr}
\Omega_r=\sum\nolimits_\xi(-1)^r\xi^{\ell-r}\omega_{\xi}
\end{equation}
($r=0,\ldots\ell$). The explicit form of the metric on $M^0$ easily follows:
we define $\pF_j(t)= \Mpc(t)\oF_j(t)$, and observe that
$\Mp'(\xi_j)=\Mpc(\xi_j)\Delta_j$, since $\Delta_j=\prod_{k\neq
j}(\xi_j-\xi_k)$.

It remains to establish the explicit form of the potentials: observe that
\begin{equation*}
\omega=\Omega_0+\sum_{r=1}^\ell d(\sigma_r \theta_r)
=dd^cH_0+\sum_{j=1}^\ell
dJ\Bigl(\frac{\xi_j^\ell\, \Mpc(\xi_j)d\xi_j}{\pF_j(\xi_j)}\Bigr)=dd^cH.
\end{equation*}
The equation $H=\sum_{r=1}^\ell \sigma_r\partial G/\partial\sigma_r-G$
determines $G$ up to a linear combination of the $\sigma_r$ with basic
coefficients. We also require $d_hd^c_hG=-\omega_h$ so the functions
$u_r=\partial G/\partial\sigma_r$ are pluriharmonic; $G$ given by~\eqref{eq:G}
has the required properties with
\begin{equation}\label{eq:ur}
u_r=-H_r-\sum_{j=1}^\ell\int^{\xi_j} \frac{(-1)^r \Mpc(t)
t^{\ell-r}}{\pF_j(t)}dt
\end{equation}
for $r=1,\ldots \ell$, where $H_r$ is a $dd^c$-potential for $\Omega_r$ on
$S$.

In the formula for $\phi$, we have used the fact that $\xi_j
\sigma_{r-1}(\hat\xi_j)=\sigma_{r}-\sigma_{r}(\hat\xi_j)$.  It is then
straightforward to check that $\Phi$ is a $dd^c$ potential for
$\phi+\sigma_1\omega$.
\end{proof}

\begin{rem}\label{r:ur} Formally, we set $u_0=-H$, so that
$G=\sum_{r=1}^\ell\sigma_ru_r-H=\sum_{r=0}^\ell\sigma_r u_r$.  Similarly, we
can write $u_{-1}=\Phi$, and in general extend~\eqref{eq:Omr}--\eqref{eq:ur}
to all $r\leq\ell$, where $dd^c H_r=\Omega_r$.  For $r>0$ $dd^c u_r=0$, while
for $r=-k\leq 0$ we have $dd^c u_{-k} = (-1)^{k+1}\phi_k$ where
\begin{align*}
\phi_k&=\sum_{\smash\xi} \xi^{\ell+k}\omega_\xi+\sum_{j,r=1}^\ell
d\Bigl(\frac{\xi_j^{\ell+k}\sigma_{r-1})(\hat\xi_j)}{\Delta_j}
Jd\xi_j\Bigr)\\
&=\sum_\xi \Bigl(\frac{\xi^{\ell+k}}{\prod_{k=1}^\ell(\xi-\xi_k)}
+\sum_{j=1}^\ell\frac{\xi_j^{\ell+k}}{\Delta_j(\xi_j-\xi)}\Bigr)
\Mpn(\xi)\omega_\xi\\
&\qquad\qquad+\sum_{i,j,r=1}^\ell\frac{\partial}{\partial\xi_i}
\Bigl(\frac{\xi_j^{\ell+k}\sigma_{r-1}(\hat\xi_j)}{\Delta_j}\Bigr)
d\xi_i\wedge\theta_r.
\end{align*}
Using~\eqref{eq:id1} and~\eqref{eq:dxi} from Appendix B, this may be written
\begin{equation*}
-J\phi_k = \sum_{s=0}^k h_{k-s}\,(-J\phi)^s
\end{equation*}
and $h_p$ is the $p$th complete symmetric function in $\xi_1,\ldots\xi_\ell$.
\end{rem}

\subsection{The hamiltonian $2$-form}

In order to complete the classification of K\"ahler metrics with a hamiltonian
$2$-form, we must show that the explicit metric of Theorem~\ref{t:one}
actually admits a hamiltonian $2$-form, with no further constraints.

\begin{thm}\label{t:two}
Let $(g,J,\omega)$ be a K\"ahler structure given explicitly
by~\eqref{eq:explicit}. Then the $J$-invariant $2$-form $\phi$ defined
by~\eqref{eq:ham?} is a hamiltonian $2$-form of order $\ell$.
\end{thm}
\begin{proof} Obviously (for the case $m=1$), the trace of $\phi$ is a Killing
potential. Hence, in order to show that $\phi$ is hamiltonian, we must show
that
\begin{equation*}
\nabla A=d\sigma_1\otimes\omega
+\tfrac12(d\sigma_1\wedge J-d^c\sigma_1\wedge\Id)
\end{equation*}
where $A=\phi+\sigma_1\omega$. Since $A$ is manifestly closed, we only need to
check the equation for $(\nabla_X A)(Y,Z)$, $(\nabla_X A)(Y,U)$, $(\nabla_U
A)(V,X)$ and $(\nabla_U A)(V,W)$, where $X,Y,Z$ and $U,V,W$ are arbitrary
horizontal and vertical vector fields respectively. Two of these equations are
immediate: $(\nabla_X A)(Y,Z)=0$ and $(\nabla_U A)(V,X)=0$.

We next consider the equation for $(\nabla_X A)(Y,U)$, which reduces,
using Proposition~\ref{p:lcc}, to the equation
\begin{equation}\label{eq:todo}
2[\phi,C(X)] = \g{X,\cdot} \otimes Jd\sigma_1^\sharp
- \omega(X,\cdot)\otimes d\sigma_1^\sharp.
\end{equation}
Here the left hand side is the commutator of $\phi$ with
\begin{equation*}
2C(X) = \sum_{r=1}^\ell \bigl(\Omega_r(X)\otimes Jd\sigma_r^\sharp
- \Omega_r(JX)\otimes d\sigma_r^\sharp\bigr).
\end{equation*}
Decomposing into the eigenspaces of $-J\phi$, we compute that
\begin{align*}
2[\phi,C(X)] &= -\sum_{r,j=1}^\ell \sum_{\smash\xi}
(-1)^r \xi^{\ell-r}\sigma_{r-1}(\hat\xi_j)
(\xi_j-\xi)
\bigl(\omega_\xi(X)\otimes d\xi_j^\sharp
+ \omega_\xi(JX)\otimes J d\xi_j^\sharp\bigr)\displaybreak[0]\\
&= -\sum_{j=1}^\ell\sum_{\smash\xi} \Bigl(\xi^{\ell}+
\sum_{r=1}^\ell(-1)^r \xi^{\ell-r}
\bigl(\xi_j\sigma_{r-1}(\hat\xi_j)+\sigma_r(\hat\xi_j)\bigr)\Bigr)\\
&\qquad\qquad\qquad\qquad\qquad\qquad\qquad\qquad
\bigl(\omega_\xi(X)\otimes d\xi_j^\sharp
+ \omega_\xi(JX)\otimes J d\xi_j^\sharp\bigr)\\
&=-\sum_{r=0}^\ell\sum_{\smash\xi} (-1)^r\sigma_r\xi^{\ell-r}
\bigl(\omega_\xi(X)\otimes d\sigma_1^\sharp
+ \omega_\xi(JX)\otimes J d\sigma_1^\sharp\bigr).
\end{align*}
This proves~\eqref{eq:todo}, by the definition of $g$ and $\omega$.

It remains to verify the equation for $(\nabla_U A)(V,W)$. Since the fibres
are totally geodesic, this amounts to showing that $\phi$ is a hamiltonian
$2$-form on the fibres. Hence it suffices to prove the result in the
orthotoric case, when we have $\phi+\sigma_1\omega=\sum_r d(f_r d^cu_r)$, with
$f_r=\sigma_1\sigma_r-\sigma_{r+1}$.

To do this we shall use only the fact that $\sum_r f_r du_r=d\Phi$ is
closed---equivalently $A$ is $J$-invariant:
$A(JK_r,K_s)=df_r(JK_s)=df_s(JK_r)=A(JK_s,K_r)$,
cf.~Proposition~\ref{p:otchar}. We define $\Phi_{rs}:=A(JK_r,K_s)$ and recall
that $\g{K_r,K_s}=H_{rs}$.

Since $A$ is closed and $J$-invariant, it suffices to check
\begin{align}\label{easy}
(\nabla_{K_r}A)(JK_s,K_t)&=0\\
\label{hard}
(\nabla_{JK_r}A)(JK_s,K_t)&=
-d\sigma_1(JK_r)\g{K_s,K_t}\\
\notag&\qquad
-\tfrac12d\sigma_1(JK_s)\g{K_r,K_t}-\tfrac12 d\sigma_1(JK_t)\g{K_s,K_r}.
\end{align}
Equation~\eqref{easy} follows immediately: the left hand side
is
\begin{equation*}
K_r\cdot\bigl(A(JK_s,K_t)\bigr)
-A(J\nabla_{K_r}K_s,K_t)-A(JK_s,\nabla_{K_r}K_t)
\end{equation*}
and all three terms are zero here, since $A(JK_s,K_t)$ is
$K_r$-invariant, $J\nabla_{K_r}K_s$ is a linear combination of the
$K_t$'s and $A(K_r,K_s)=0$ for all $r,s$.

On the other hand, for equation~\eqref{hard} we have
\begin{equation*}
(\nabla_{JK_r}A)(JK_s,K_t)=
JK_r\cdot\bigl(df_s(JK_t)\bigr)+df_t(\nabla_{K_r}K_s)+df_s(\nabla_{K_r}K_t).
\end{equation*}
Now
\begin{equation*}
\nabla_{K_r}K_s=-\frac12\grad\g{K_r,K_s}=-\frac12\sum_p\frac{\partial
H_{rs}} {\partial u_p}du_p^\sharp
\end{equation*}
and $H_{rs}=\partial^2H/\partial u_r\partial u_s$, so~\eqref{hard}
holds if and only if
\begin{equation}\label{rtp}
\frac{\partial\Phi_{st}}{\partial u_r}-\frac12\sum_{p,q}G_{pq}
\Bigl(\Phi_{tp}\frac{\partial H_{sq}}{\partial u_r}+\Phi_{sp}
\frac{\partial H_{tq}}{\partial u_r}\Bigr)
=H_{1r}H_{st}+\frac12(H_{1s}H_{rt}+H_{1t}H_{rs}).
\end{equation}
This simplifies once we observe that $\Phi_{st}
=\sum_{p,q} G_{pq}\Phi_{tp}H_{sq}=\sum_{p,q} G_{pq}\Phi_{sp}H_{tq}$,
so that the product rule reduces the left hand side to
\begin{equation}\label{lhs}
\frac12\sum_{p,q}\Bigl(H_{sq}\frac{\partial}{\partial u_r}(G_{pq}\Phi_{tp})
+H_{tq}\frac{\partial}{\partial u_r}(G_{pq}\Phi_{sp})\Bigr).
\end{equation}
Now we use the fact that $f_s=\sigma_1\sigma_s-\sigma_{s+1}$ and
$\partial\sigma_s/\partial u_r=H_{rs}$ to deduce that
\begin{equation*}
\sum_p\frac{\partial}{\partial u_r}(G_{pq}\Phi_{tp})=
H_{rt}\delta_{1q}+H_{1r}\delta_{tq}.
\end{equation*}
Substituting this into~\eqref{lhs} yields~\eqref{rtp}, and hence~\eqref{hard}.
\end{proof}

\section{The curvature of K\"ahler manifolds with a hamiltonian $2$-form}
\label{s:ckmh}

\subsection{The Ricci potential and scalar curvature}

In this section we compute the Ricci potential and scalar curvature for any
K\"ahler manifold $(M,g,J,\omega)$ given by~\eqref{eq:explicit}, using the
formulae obtained in section~\ref{s:hta}.  In terms of the Vandermonde matrix
$V_{rj}$ and its inverse $W_{jr}$ (see Appendix B), we have
\begin{equation}
G_{rs}=\sum_{j=1}^\ell
\frac{\Mpc(\xi_j)V_{rj}V_{sj}}{\pF_j(\xi_j)\Delta_j},\qquad
H_{rs}=\sum_{j=1}^\ell \frac{W_{jr}W_{js}\pF_j(\xi_j)\Delta_j}{\Mpc(\xi_j)},
\end{equation}
and hence, up to a sign, $\det G_{rs}$ is $\prod_{j=1}^\ell\Mpc(\xi_j)
\pF_j(\xi_j)^{-1}$.  Also, the formula $\omega_h=
\sum_{\xi}\Mpn(\xi)\omega_\xi$ for the K\"ahler quotient gives $Q=\prod_\xi
\Mpn(\xi)^{m_\xi} =\pm\prod_{j=1}^{\smash\ell}\Mpc(\xi_j)$,
cf.~\eqref{eq:Q}. It follows immediately from Proposition~\ref{p:pluriricci}
that if $\kappa_\xi$ is a Ricci potential for $(S_\xi,\omega_\xi)$, then a
Ricci potential for $(M,\omega)$ is
\begin{equation}\label{eq:rp}
\kappa=\sum_\xi\kappa_\xi -\frac12\sum_{j=1}^\ell \log |\pF_j(\xi_j)|.
\end{equation}

In order to obtain the scalar curvature from this, we need a formula for the
laplacian in the $\xi_j$ coordinates.

\begin{lemma} For any function $f$, we have
\begin{equation}\label{eq:lapfxi}
\Delta f = \Delta_h f -\sum_{j=1}^\ell \frac 1{\Delta_j\,\Mpc(\xi_j)}
\frac{\partial}{\partial\xi_j}\Bigl( \pF_j(\xi_j) \frac{\partial f}
{\partial \xi_j}\Bigr).
\end{equation}
\end{lemma}
\begin{proof} We just need to change coordinates in equation~\eqref{eq:lapf}.
For this observe that
\begin{gather*}\notag
\frac{\partial}{\partial\sigma_r}=\sum_{k=1}^\ell
\frac{V_{rk}}{\Delta_k}\frac{\partial}{\partial\xi_k}
\qquad\text{and}\qquad\frac{\partial}{\partial u_r}
=\sum_{s=1}^\ell H_{rs} \frac{\partial}{\partial\sigma_s}
= \sum_{j=1}^\ell \frac{W_{jr}\pF_j(\xi_j)}{\Mpc(\xi_j)}
\frac{\partial}{\partial\xi_j},\\
\tag*{so that}
\Delta f = \Delta_h f - \frac1{Q} \sum_{r,j,k=1}^\ell
\frac{V_{rk}}{\Delta_k}\frac{\partial}{\partial\xi_k}
\Bigl(\frac{Q W_{jr} \pF_j(\xi_j)}{\Mpc(\xi_j)}\,
\frac{\partial f}{\partial\xi_j}\Bigr).
\end{gather*}
Since $Q=\pm \prod_{i=1}^\ell \Mpc(\xi_i)$, this agrees with~\eqref{eq:lapfxi}
once we observe that
\begin{equation*}
\sum_{k,r=1}^\ell
\frac{V_{rk}}{\Delta_k}\frac{\partial W_{jr}}{\partial\xi_k}
=-\sum_{r,k=1}^\ell
\frac{W_{jr}}{\Delta_k}\frac{\partial V_{rk}}{\partial\xi_k}
=\sum_{r=1}^{\ell-1} \sum_{k=1}^\ell\frac{W_{jr}}
{\Delta_k}(-1)^r(\ell-r)\xi_k^{\ell-r-1}=0
\end{equation*}
by the Vandermonde identity.
\end{proof}
Applying this formula to the Ricci potential, we deduce immediately that
\begin{equation}\label{eq:explicitscal}
\Scal=\sum_{\smash \xi} \frac{\Scal_{g_\xi}}{\Mpn(\xi)}
-\sum_{j=1}^\ell \frac{\pF_j''(\xi_j)}{\Delta_j\, \Mpc(\xi_j)},
\end{equation}
where $\Scal_{g_\xi}$ is the scalar curvature of the (possibly negative
definite) metric $g_\xi$.

\begin{lemma}\label{l:scal} Suppose that $\Scal$ depends polynomially on
$\xi_1,\ldots\xi_\ell$. Then
\begin{itemize}
\item for all $j$, $\pF_{j}''(t)=\check \Mpc(t)R(t)$, where $\check
\Mpc(t)=\prod_{\xi} (t-\xi)^{m_\xi-1}_{\smash[b]{\phantom{T}}}$ and $R(t)$ is
a polynomial independent of $j$\textup;
\item for all $\xi$, $(g_\xi,\omega_\xi)$ has
$\Scal_{g_\xi}=-R(\xi)/\prod_{\eta\neq\xi} (\xi-\eta)$.
\end{itemize}
We then have
\begin{equation}\label{eq:nicescal}
\Scal=-\sum_\xi \frac{R(\xi)}{\prod_{\eta\neq\xi}(\xi-\eta)
\prod_k (\xi-\xi_k)}-\sum_{j=1}^\ell \frac{R(\xi_j)}{\prod_{k\neq j}
(\xi_j-\xi_k)\prod_\eta(\xi_j-\eta)},
\end{equation}
where $\check m=\ell+\sum_\xi 1 = m - \sum_\xi (m_\xi-1)$.

Furthermore, if $\Scal$ has degree $\leq q$ in each variable $\xi_j$, then
$R(t)$ has degree at most $\check m+q-1$. Hence, for each $j$, $\pF_j(t)$ is a
polynomial of degree at most $m+q+1$.
\end{lemma}
\begin{proof}
We multiply the formula~\eqref{eq:explicitscal} by $\Delta_k\,\Mpc(\xi_k)$
to obtain an equality between polynomials in $\xi_k$ (on a nonempty open
set, hence everywhere):
\begin{multline*}\notag
\Delta_k\,\Mpc(\xi_k)\,\Scal\\
=-\sum_\xi \frac{\Delta_k\, \check \Mpc(\xi_k)\prod_{\eta\neq\xi}(\xi_k-\eta)}
{\prod_{j\neq k}(\xi-\xi_j)}\Scal_{g_\xi}
-\pF_k''(\xi_k)-\sum_{j\neq k} \frac{\Delta_k \,\Mpc(\xi_k)}
{\Delta_j \,\Mpc(\xi_j)} \pF_j''(\xi_j).
\end{multline*}
This clearly shows that $\pF_k''$ is a polynomial with $\check \Mpc$ as a
factor. Evaluating at $\xi_k=\xi_j$ for some fixed $j$, we obtain
$\pF_k''(\xi_j)=\pF_j''(\xi_j)$ for all $\xi_j$ (in a nonempty open set, hence
everywhere). Dividing through by $\check \Mpc(\xi_k)$ we now have
\begin{multline*}\notag
\Delta_k\prod_{\eta}(\xi_k-\eta)\,\Scal\\
=-\sum_\xi \frac{\Delta_k\prod_{\eta\neq\xi}(\xi_k-\eta)}
{\prod_{j\neq k}(\xi-\xi_j)}\Scal_{g_\xi}-R(\xi_k)-\sum_{j\neq k}\frac{\Delta_k
\prod_\eta(\xi_k-\eta)}{\Delta_j\prod_\eta (\xi_j-\eta)}R(\xi_j).
\end{multline*}
Evaluating at $\xi_k=\xi$ gives the formula for $\Scal_{\smash{g_\xi}}$, and
it is straightforward to count the degree in $\xi_k$. Dividing
by $\Delta_k\,\prod_{\eta}(\xi_k-\eta)$ now gives~\eqref{eq:nicescal}.
\end{proof}
To interpret~\eqref{eq:nicescal}, we adjoin the distinct constant roots to
the variables $\xi_1,\ldots\xi_\ell$.  If we label these $\xi_1,\ldots
\xi_\ell,\xi_{\ell+1},\ldots \xi_{\check m}$, and let
$\Delta_j^\vee=\prod_{k\neq j}(\xi_j-\xi_k)$, where the product is over
$k=1,\ldots \check m$, then the right hand side of~\eqref{eq:nicescal} is just
$-\sum_{j=1}^{\smash{\check m}} R(\xi_j)/\Delta_j^\vee$ which is a polynomial
of degree at most $q$ in each $\xi_j$, by the Vandermonde identity.

\subsection{Extremal K\"ahler metrics}

Recall that a K\"ahler metric is called {\it extremal} if the scalar curvature
is a Killing potential~\cite{calabi1}.  Weakly Bochner-flat K\"ahler metrics
of dimension $2m\geq4$ are extremal, since the scalar curvature is then the
trace of a hamiltonian $2$-form. In this section we classify the extremal
K\"ahler metrics with a hamiltonian $2$-form such that $d_h\Scal=0$.

\begin{prop}\label{p:ext}
Let $(M,g,J,\omega)$ be K\"ahler with a hamiltonian $2$-form $\phi$. Then
$\Scal$ is a hamiltonian for a Killing vector field tangent to the fibres of
$M$ over the K\"ahler quotient $S$ if and only if $(g,J,\omega)$ has the
explicit form~\eqref{eq:explicit} where\textup:
\begin{itemize}
\item for all $j$, $\pF_{j}''(t)=\check \Mpc(t)\bigl( \sum_{r=0}^{\check
m} a_r t^{\check m-r}\bigr)$, and $a_0,\ldots a_{\check m}$ are arbitrary
constants \textup(independent of $j$\textup)\textup;
\item for all $\xi$, $(g_\xi,\omega_\xi)$ has
$\Scal_{\smash{g_\xi}}=-\bigl(\sum_{r=0}^{\check m} a_r \xi^{\check m-r}\bigr)
/\prod_{\eta\neq\xi} (\xi-\eta)$.
\end{itemize}
The scalar curvature of $(g,J,\omega)$ is then given by
$\Scal=-(a_0\check\sigma_1+a_1)$, where $\check\sigma_1:=\sum_{j=1}^{\check
m}\xi_j=\sigma_1+\sum_\xi\xi$, so that $\Scal$ is a hamiltonian for $-a_0
K_1$.

Any constant scalar curvature K\"ahler metric with a hamiltonian $2$-form
arises in this way with $a_0=0$, and is scalar-flat if and only if $a_1=0$.
\end{prop}
\begin{proof} Since $\Scal$ is invariant under $K_1,\ldots K_\ell$, it must
be a function of $\sigma_1,\ldots \sigma_\ell$, and since $J\grad_g \Scal$
commutes with $K_1,\ldots K_\ell$ and is in their span at each point, it must
in fact be a constant linear combination of $K_1,\ldots K_\ell$, so that
$\Scal$ is an affine function of $\sigma_1,\ldots \sigma_\ell$. Now any such
function is a polynomial in $\xi_1,\ldots \xi_\ell$ of degree one in each
$\xi_j$.  Hence we can apply Lemma~\ref{l:scal}.

Conversely, the Vandermonde identities imply that
$\Scal=-(a_0\check\sigma_1+a_1)$, which is a Killing potential for $-a_0K_1$,
since $\check\sigma_1$ differs from $\sigma_1$ by a constant.
\end{proof}

\subsection{Weakly Bochner-flat K\"ahler metrics}

On a weakly Bochner-flat K\"ahler manifold of dimension $2m\geq4$, $\nrho$
is a hamiltonian $2$-form, so we obtain a classification by specializing the
work of the previous section to the case $\nrho=\phi$.  In fact we may as well
consider more generally the case that $\nrho=a\phi+b\omega$ for constants
$a,b$. Then when $a=0$, we will have characterized the K\"ahler--Einstein
metrics among K\"ahler metrics with a hamiltonian $2$-form. Note however, that
we have fixed $\phi$: it is obviously possible for $\nrho$ to be a hamiltonian
$2$-form without it being equal to $a\phi+b\omega$ for a given $\phi$, but we
have nothing to say about this situation.

\begin{prop}\label{p:wbf}
Let $(M,g,J,\omega)$ be K\"ahler of dimension $2m\geq4$ with a hamiltonian
$2$-form $\phi$. Then $M$ is weakly Bochner-flat with $\nrho$ a linear
combination of $\phi$ and $\omega$ if and only if $(g,J,\omega)$ has the
explicit form~\eqref{eq:explicit} where\textup:
\begin{itemize}
\item for all $j$, $\pF_j'(t)= \Mpc(t)\bigl(\sum_{r=-1}^{\ell} b_r
t^{\ell-r}\bigr)$, and $b_{-1},\ldots b_{\ell}$ are arbitrary constants
\textup(independent of $j$\textup)\textup;
\item for all $\xi$, $(g_\xi,\omega_\xi)$ is K\"ahler--Einstein with
K\"ahler--Einstein constant
\begin{equation*}
\frac1{m_\xi}\Scal_{g_\xi}=-\sum_{r=-1}^\ell b_r\xi^{\ell-r}.
\end{equation*}
\end{itemize}
The Ricci form of $(g,J,\omega)$ is then given by
$\rho=-\frac12\bigl(b_{-1}(\phi+\sigma_1\omega)+b_0\omega\bigr)$.

Any K\"ahler--Einstein metric with a hamiltonian $2$-form arises in this way
with $b_{-1}=0$, and is Ricci-flat if and only if $b_0=0$.
\end{prop}
\begin{proof}
$-2\rho=b_{-1}(\phi+\sigma_1\omega)+b_0\omega$ if and only if $dd^c$
potentials for the two sides differ by a pluriharmonic function. This means
that $-2\kappa$ must be of the form $-\sum_{r=-1}^\ell (-1)^rb_ru_r$ where
\begin{equation*}
u_r=-H_r-\sum_{j=1}^\ell\int^{\xi_j} \frac{(-1)^r t^{\ell-r}\Mpc(t)}
{\pF_j(t)}dt
\end{equation*}
for $r=-1,\ldots \ell$, the $H_r$ being functions on $\prod_\xi S_\xi$ such
that $H_r$ is a $dd^c$-potential for $\Omega_r=\sum_\xi
(-1)^r\xi^{\ell-r}\omega_\xi$ when $r=0,\ldots\ell$. (So $u_{-1}$ is a $dd^c$
potential for $\phi+\sigma_1\omega$, $u_0$ is a $dd^c$ potential for $-\omega$
and $u_r$ is pluriharmonic for $r\geq 1$---see Remark~\ref{r:ur}.) Now it
follows from~\eqref{eq:rp} that
\begin{equation*}
-2\kappa=-2\sum_\xi\kappa_\xi + \sum_{j=1}^\ell \int^{\xi_j}
\frac{\pF_j'(t)}{\pF_j(t)} dt.
\end{equation*}
Comparing this with
\begin{equation*}
\sum_{r=-1}^\ell (-1)^r b_r H_r +
\sum_{j=1}^\ell\int^{\xi_j} \frac{\Mpc(t)\bigl(
\sum_{r=-1}^\ell b_r t^{\ell-r}\bigr)}{\pF_j(t)}dt,
\end{equation*}
we obtain the required result.
\end{proof}

\subsection{Bochner-flat K\"ahler metrics}\label{s:bf}

We now rederive Bryant's classification of Bochner-flat K\"ahler metrics in
the present framework. One interesting feature of our approach is that we can
at the same time explicitly classify hamiltonian $2$-forms on Bochner-flat
K\"ahler--Einstein manifolds, i.e., on K\"ahler manifolds of constant
holomorphic sectional curvature, cf.~section~\ref{s:chsc}.

We set $\hat \Mpc(t)=\prod_\xi (t-\xi)^{m_\xi+1}$ and $\hat m=\ell-\sum_\xi
1=m-\sum_\xi (m_\xi+1)$.

\begin{prop}\textup{\cite{bryant}}\label{p:bf} Let $(M,g,J,\omega)$ be a
K\"ahler manifold of dimension $2m\geq4$ with a hamiltonian $2$-form.  Then
$M$ is Bochner-flat with $\nrho$ a linear combination of $\phi$ and $\omega$
if and only if $(g,J,\omega)$ has the explicit form~\eqref{eq:explicit}
where\textup:
\begin{itemize}
\item for all $j$, $\pF_j(t)= \hat \Mpc(t)\bigl(\sum_{r=-2}^{\hat m} c_r
t^{\hat m-r}\bigr)$, and $c_{-2},\ldots c_{\hat m}$ are arbitrary constants
\textup(independent of $j$\textup)\textup;
\item for all $\xi$, $(g_\xi,\omega_\xi)$ has constant holomorphic sectional
curvature
\begin{equation*}
\frac1{m_\xi(m_\xi+1)}\Scal_{g_\xi}
=-\Bigl(\sum_{r=-2}^{\hat m} c_r\xi^{\hat m-r}\Bigr)
\prod_{\eta\neq\xi}(\xi-\eta).
\end{equation*}
\end{itemize}
The curvature of $(g,J,\omega)$ is then given by $R=-\{J\circ\hat\rho,\cdot\}
+\hat\rho\otimes\omega+\omega\otimes\hat\rho$, where
$\hat\rho=-\frac12\bigl(c_{-2}(\phi+\frac12\hat\sigma_1\omega)+\frac12
c_{-1}\omega\bigr)$ and $\hat\sigma_1=\sigma_1-\sum_\xi\xi$.

Any constant holomorphic sectional curvature K\"ahler metric with a
hamiltonian $2$-form arises in this way with $c_{-2}=0$, and is flat if and
only if $c_{-1}=0$.
\end{prop}
\begin{proof} By Proposition~\ref{p:wbf}, we may assume that $(g,J,\omega)$
is weakly Bochner-flat, with $\pF_j'$ a polynomial divisible by $\Mpc$ and
$S_\xi$ K\"ahler--Einstein. Since $\nrho$ is a linear combination of $\phi$
and $\omega$, we may use the system~\eqref{nablaphi2}.  Integrating the last
three equations gives $\tau_0=C_{-2}$, $\tau_1= C_{-2}\sigma_1+C_{-1}$ and
$\tau_2=C_{-2}(\sigma_1^2-\sigma_2)+C_{-1}\sigma_1+C_0$ for constants
$C_{-2}$, $C_{-1}$ and $C_0$. (We have used the fact that $\sigma-\sigma_1$
and $\g{\phi,\phi}+2\sigma_2-\sigma_1^2$ are constants.) Hence
\begin{align*}\notag
dd^c\sigma_1 &= \frac2m \WK(\phi)
-C_{-2}J\phi^2+(C_{-2}\sigma_1+C_{-1})\phi+(C_{-2}(\sigma_1^2-\sigma_2)
+C_{-1}\sigma_1+C_0)\omega\\
&=\frac2m \WK(\phi)-C_{-2}(J\phi^2-\sigma_1\phi-(\sigma_1^2-\sigma_2)\omega)
+C_{-1}(\phi+\sigma_1\omega)+C_0\omega.
\end{align*}
It follows from Remark~\ref{r:ur} that
\begin{equation*}
\frac2m \WK(\phi)=dd^c(\sigma_1+C_{-2}u_{-2}-C_{-1}u_{-1}+C_0u_0).
\end{equation*}
(Here we recall that $dd^c u_{-2}= J \phi^2-\sigma_1 \phi - (\sigma_1^2 -
\sigma_2) \omega$.)

Thus $\WK(\phi)$ is basic if and only if $\pF_j(t)=\Mpc(t)\bigl(
\sum_{r=-2}^\ell C_r t^{\ell-r}\bigr)$, where $C_1,\ldots C_\ell$ are
arbitrary constants. This means that $\pF_j$ is a polynomial divisible by
$\Mpc$. Since $\pF_{\smash j}'$ is also divisible by $\Mpc$, $\WK(\phi)$ is
basic if and only if $\pF_j$ is a polynomial divisible by $\hat \Mpc$. In
particular $\sum_{r=-2}^{\smash\ell} C_r \xi^{\ell-r}=0$ for each constant
root $\xi$, and so
\begin{align*}
\frac2m \WK(\phi)&=dd^c\Bigl(\sigma_1+\sum_{r=-2}^\ell (-1)^r
C_r u_r\Bigr)\\
&=-\sum_{r=-2}^\ell (-1)^r C_r\, dd^cH_r = \sum_{r=-2}^\ell \sum_\xi C_r
\xi^{\ell-r} \omega_\xi  =0.
\end{align*}
Hence we may suppose that $\pF_j(t)=\hat\Mpc(t)\bigl(\sum_{r=-2}^{\hat m} c_r
t^{\hat m-r}\bigr) $ and that $\WK(\phi)=0$. Note that $c_{-2}=C_{-2}$ and
$c_{-1}=C_{-1}+(\sum_\xi\xi)C_{-2}$, so $\tau_0= c_{-2}$ and
$\tau_1=c_{-2}\hat\sigma_1+c_{-1}$.

Since $\nrho$ is an constant linear combination of $\phi$ and $\omega$,
equation~\eqref{curvphi} implies that $[\WK(\psi),\phi]= \frac 1m
[\WK(\phi),\psi]=0$. Also equation~\eqref{eq:WKK1} implies that $\iota _{K}
\WK(\psi)=0$.  Now any $2$-form commuting with $\phi$ is the sum of a vertical
and a horizontal $2$-form (i.e., there is no mixed component); then since
$-J\phi$ has distinct eigenvalues on the fibres of $M$ over $S$, the vertical
component is of the form $\sum_{j=1}^{\ell} \mu _j d \xi_j \wedge J d \xi
_j$. If also the contraction with $K$ is zero, we have $\sum_{j=1}^{\smash
\ell} \mu _j d \xi _j (X) |d \xi _j| ^2 = 0$, for each vector field $X$, which
forces $\mu _1 = \cdots = \mu _r = 0$. We deduce that for any $2$-form $\psi$,
$\WK(\psi)$ is horizontal, and, since $\WK$ is symmetric, $\WK(\psi)$ vanishes
unless $\psi$ is horizontal.

We next employ the Gray--O'Neill submersion formulae~\cite{Gray,ONeill} which
apply in this situation (the submersion of $M$ over $S$ is not riemannian, so
we need the framework of Gray). If $X,Y,Z$ are horizontal vector fields, we
have (see Proposition~\ref{p:lcc})
\begin{align*}
(R_{X,Y} Z)^H &= R^H_{X,Y} Z -\g{C(X),C(Y,Z)}^\sharp+\g{C(Y),C(X,Z)}^\sharp\\
 &\qquad\qquad\;\;+\g{C(Z),C(X,Y)}^\sharp-\g{C(Z),C(Y,X)}^\sharp,\\
\tag*{with}
2C(X,Y) &= \sum_{r=1}^\ell \bigl(\Omega_r(X,Y) K_r + \Omega_r(JX,Y) JK_r\bigr).
\end{align*}
To compute $\g{C(X,\tilde X), C(Y,\tilde Y)}$ for horizontal vector
fields $X,\tilde X,Y,\tilde Y$, we use the definition of $C$, and expand
$\Omega_r$ and $\g{K_r,K_s}$ to get
\begin{equation*}\notag
\g{C(X,\tilde X),C(Y,\tilde Y)} = \frac14\sum\nolimits_{\xi,\eta}
f_{\xi,\eta}\, \bigl(\omega_\xi(X,\tilde X)\omega_\eta(Y,\tilde Y)+
\omega_\xi(JX,\tilde X)\omega_\eta(JY,\tilde Y)\bigr),
\end{equation*}
where
\begin{align*}
f_{\xi,\eta} &= \sum_{j=1}^\ell \frac{\pF_j(\xi_j)}{\Mpc(\xi_j)\Delta_j}
\Bigl(\sum_{r=1}^\ell (-1)^r\xi^{\ell-r}\sigma_{r-1}(\hat\xi_j)\Bigr)
\Bigl(\sum_{s=1}^\ell (-1)^s\eta^{\ell-s}\sigma_{s-1}(\hat\xi_j)\Bigr)\\
&=\Bigl(\prod_{k=1}^\ell (\xi-\xi_k)(\eta-\xi_k)\Bigr)\Bigl(
\sum_{j=1}^\ell \frac{\pF_j(\xi_j)}
{\Mpc(\xi_j)\Delta_j(\xi_j-\xi)(\xi_j-\eta)}\Bigr)\\
&=\Mpn(\xi)\Mpn(\eta)\biggl(
C_{-2}(\sigma_1+\xi+\eta)+C_{-1}-\delta_{\xi,\eta}
\frac{\bigl(\sum_{r=-2}^{\hat m} c_r \xi^{\hat m-r}\bigr)
\prod_{\eta\neq\xi}(\xi-\eta)}{\prod_{j=1}^\ell(\xi-\xi_j)}\biggr).
\end{align*}
For the last line of this calculation, we observe that $\pF_j(t)/\Mpc(t)$ is a
polynomial of degree $\ell+2$ vanishing when $t=\xi$; then if $\xi\neq\eta$ we
can apply the Vandermonde identity with the variables
$\xi_1,\ldots\xi_\ell,\xi,\eta$, whereas if $\xi=\eta$, we may use the
Vandermonde identity for the polynomial $\pF_j(t)/\bigl(\Mpc(t)(t-\xi)\bigr) =
\bigl(\sum_{r=-2}^{\hat m} c_r t^{\hat m-r}\bigr) \prod_{\eta\neq
\xi}(t-\eta)$ of degree $\ell+1$, with the variables
$\xi_1,\ldots\xi_\ell,\xi$.

We expand the curvature $R$ using~\eqref{eq:newdecRK} and the fact that
\begin{equation*}\notag
\hat\rho = -\tfrac12(\tau_0\phi+\tfrac12\tau_1\omega)
=-\tfrac12 C_{-2}\bigl(\phi+\tfrac12\sigma_1\omega\bigr)-\tfrac14C_{-1}\omega
=-\tfrac12 c_{-2}\bigl(\phi+\tfrac12\hat\sigma_1\omega\bigr)
-\tfrac14c_{-1}\omega.
\end{equation*}
(See section~\ref{s:dswbf}.) The final ingredient in the computation, from
Proposition~\ref{p:wbf}, is the fact that each $S_\xi$ is K\"ahler--Einstein
with K\"ahler--Einstein constant
\begin{equation*}
\frac1{m_\xi}\Scal_{g_\xi} = -(m_\xi+1)
\Bigl(\sum_{r=-2}^{\smash{\hat m}} c_r\xi^{\hat m-r}\Bigr)
\prod_{\eta\neq\xi}(\xi-\eta).
\end{equation*}
Putting these ingredients together, bearing in mind that
$\smash{\omega\restr{S_\xi}}=\Mpn(\xi)\omega_\xi$, we find that a remarkable
cancellation occurs (cf.~\cite{bryant}) and we obtain
\begin{equation*}
(\WK_{X,Y}Z)^H = \sum\nolimits_{\smash{\xi}} W^{\mathcal K,S_\xi}_{X,Y} Z,
\end{equation*}
where $W^{\mathcal K,S_\xi}$ denotes the Bochner tensor of $S_\xi$ (pulled
back to the K\"ahler product). We deduce that $(\WK_{X,Y}Z)^H=0$ if and only
if each $S_\xi$ has constant holomorphic sectional curvature given by the
stated formula.
\end{proof}

Our proof above is very much inspired by~\cite[Section 4.5]{bryant}, where
Bryant indicates how to obtain an explicit formula for the general
Bochner-flat metric, although he stops short of providing the final
formula. Our approach has proceeded in reverse, by first finding the general
formula, then showing that it is Bochner-flat. This has permitted us to give a
proof using standard methods in K\"ahler geometry, substituting a linear
system for the nonlinear system which Bryant integrates using Cartan's
generalization of Lie's Third Theorem.

Bryant's remarkable paper also addresses global questions: the compact
Bochner-flat K\"ahler manifolds are necessarily locally symmetric, but Bryant
finds compact orbifold examples, and classifies the complete examples
(cf.~also~\cite{AG2} for the case of K\"ahler surfaces). In forthcoming work,
we shall find that there are many compact weakly Bochner-flat K\"ahler
manifolds (see also~\cite{ACG0}).

We now return to the characteristic polynomial of section~\ref{s:dswbf}. We
write $\pF$ for $\pF_j$ (which is independent of $j$), and define a {\it
minimal polynomial} $\Fm(t):=\pF(t)/\Mpc(t)= \bigl(\sum_{r=-2}^{\hat m} c_r
t^{\hat m-r}\bigr)\prod_\xi (t-\xi)$.  When $c_{-2}\neq0$, these polynomials
are (up to affine transformation of $t$) Bryant's characteristic and reduced
characteristic polynomials.

\begin{prop} Let $\phi$ be a hamiltonian $2$-form on a Bochner-flat K\"ahler
manifold $(M,g,J,\omega)$ with $\nrho$ a linear combination of $\phi$ and
$\omega$, as in Proposition~\textup{\ref{p:bf}}. Then $\pF$ is the
characteristic polynomial $\CP$ of $(g,J,\omega,\phi)$.
\end{prop}
\begin{proof} To compute $(\tau_0t^2+\tau_1t+\tau_2)\Mp(t)-\g{K,K(t)}$,
observe that
\begin{equation*}\notag
\g{K,K(t)}=-\sum_{j=1}^\ell\Mpc(t)\Bigl(\prod_{k\neq j}(t-\xi_j)\Bigr)
|d\xi_j|^2
= \Mp(t)\sum_{j=1}^\ell \frac{\Fm(\xi_j)}{(\xi_j-t)\Delta_j}.
\end{equation*}
Hence, writing $\Fm(t)=\sum_{r=-2}^\ell C_r t^{\ell-r}$, Vandermonde
identities (with variables $\xi_1,\ldots\xi_\ell,t$) give
\begin{multline*}\notag
\g{K,K(t)}=\Mp(t)\bigl(C_{-2}((\sigma_1+t)^2-\sigma_2-t\sigma_1)
+ C_{-1}(\sigma_1+t)+C_0\bigr) - \Mp(t)\frac{\smash{\Fm(t)}}{\Mpn(t)}\\
=\Mp(t)\bigl( C_{-2}t^2 + (C_{-2}\sigma_1+C_{-1})t +
(C_{-2} (\sigma_1^2-\sigma_2) + C_{-1}\sigma_1+C_0) - \pF(t).
\end{multline*}
Hence $(\tau_0t^2+\tau_1t+\tau_2)\Mp(t)-\g{K,K(t)}=\pF(t)$.
\end{proof}

When $c_{-2}=0$ (i.e., $\tau_0=0$), Proposition~\ref{p:bf} provides a
classification of hamiltonian $2$-forms on simply-connected manifolds of
constant holomorphic sectional curvature in terms of two polynomials $\Mpc$
and $\Fm$ respectively of degrees $m-\ell$ (precisely) and $\ell+1$ (at most),
such that every root of $\Mpc$ is a root of $\Fm$.

In section~\ref{s:chsc} we showed that hamiltonian $2$-forms are then given by
parallel sections of a flat connection on a bundle of rank $(m+1)^2$. In the
simply-connected case, or when the K\"ahler manifold is an open subset of
$\C^m$, $\C P^m$ or $\C\mathcal{H}^m$, the bundle is trivial, and parallel
sections extend globally to $\C^m$, $\C P^m$ or $\C\mathcal{H}^m$. For
$s\neq0$, we identified the solution space with the Lie algebra
$\mathfrak{u}(m+1)$ or $\mathfrak{u}(m,1)$ of Killing potentials, and we gave
an explicit description for $s=0$.

The positive-definiteness of the explicit metric $g$ implies that
$\Fm(\xi_j)/\Delta_j$ must be positive for all $j$, so $\Fm$ must have at
least $\ell-1$ distinct roots (without loss of generality $\xi_1<\cdots
<\xi_\ell$ and $\Fm$ has a sign change in each interval).  The scalar
curvature is nonzero when $\Fm$ has exactly degree $\ell+1$. In the positive
case ($C_{-1}<0$), it changes sign in $(-\infty,\xi_1)$ and
$(\xi_\ell,\infty)$ and so has $\ell+1$ distinct roots; the roots of $\Mpc$
are thus the multiple roots of the characteristic polynomial $\pF$, and the
order $\ell$ of $\phi$ equals the number of different roots of the
characteristic polynomial minus one. (The case of only one root is the trivial
case that $\tau$ is constant and $\phi=\tau\omega$.)  The negative case is
more complicated, corresponding to the fact that $\mathfrak{u}(m,1)$ has
non-semisimple elements, and $\pF$ alone is not enough to classify them---we
also need to know the factorization $\pF(t)=\Fm(t)\Mpc(t)$, i.e., the minimal
polynomial.

\subsection{The Calabi-type case}

To illustrate the conditions of the previous subsections, and for use in the
next subsection, we specialize to the simple, but important case of
hamiltonian $2$-forms of order one, when the K\"ahler structure may be
written:
\begin{equation*}\begin{split}
g&={\textstyle\sum_\xi} (z-\xi)g_\xi +\frac{\prod_\xi(z-\xi)^{m_\xi}
}{\pF(z)}\, dz^2+
\frac{\pF(z)}{\prod_\xi(z-\xi)^{m_\xi}}\,\theta^2,\\
\omega&= {\textstyle\sum_\xi} (z-\xi)\omega_\xi
+dz\wedge\theta,\qquad\qquad d\theta={\textstyle\sum_\xi}\omega_\xi,\\
Jdz &= \frac{\pF(z)}{\prod_\xi(z-\xi)^{m_\xi}}\,\theta,
\qquad\qquad\qquad\quad
J\theta=-\frac{\prod_\xi(z-\xi)^{m_\xi}}{\pF(z)}\,dz,
\end{split}\end{equation*}
where $2m_\xi=\dim S_\xi$ and $\pm(g_\xi,\omega_\xi)$ is a K\"ahler metric on
$S_\xi$ (compared to~\eqref{eq:explicit}, we have reversed the sign of
$(g_\xi,\omega_\xi)$). Suppose there are $N$ different constant roots
$\xi$. Then Propositions~\ref{p:ext},~\ref{p:wbf} and~\ref{p:bf} specialize as
follows:
\begin{numlist}
\item $g$ is extremal when $P_2(t):=\pF''(t)/\prod_\xi(t-\xi)^{m_\xi-1}$ is a
polynomial of degree $N+1$ and $g_\xi$ has constant scalar curvature
$P_2(\xi)/\prod_{\eta\neq\xi}(\xi-\eta)$;

\indent\indent $g$ has constant scalar curvature when $P_2$ has degree $N$;

\item $g$ is weakly Bochner-flat when
$P_1(t):=\pF'(t)/\prod_\xi(t-\xi)^{m_\xi}$ is a polynomial of degree $2$ and
$g_\xi$ is K\"ahler--Einstein with K\"ahler--Einstein constant $P_1(\xi)$;

\indent\indent $g$ is K\"ahler--Einstein when $P_1$ has degree $1$;

\item $g$ is Bochner-flat when $P_0(t):=\pF(t)/\prod_\xi(t-\xi)^{m_\xi+1}$ is
a polynomial of degree $3-N$ and $S_\xi$ has constant holomorphic sectional
curvature $P_0(\xi)\prod_{\eta\neq\xi}(\eta-\xi)$;

\indent\indent $g$ has constant holomorphic sectional curvature
when $P_0$ has degree $2-N$.
\end{numlist}
Metrics of this form have been used in many places to provide compact or
complete examples on (projective) line bundles over K\"ahler products:
see~\cite{koisak,wang} for K\"ahler--Einstein metrics,
and~\cite{calabi1,hwang-singer,christina,christina2} for extremal
K\"ahler metrics. We shall explore and generalize these global examples in
subsequent work.

The case $N=1$ is particularly important: the constant roots are all equal, so
without loss of generality we may take them all to be zero. The polynomials
$P_0,P_1,P_2$ are all quadratic: in case (i), we may write $\pF(t)=
c_{-2}t^{m+2}+c_{-1}t^{m+1}+c_0t^m+bt+a$ and the scalar curvature of $g_0$ is
$m(m-1)c_0$; then case (ii) has $b=0$ and K\"ahler--Einstein constant $m c_0$,
while case (iii) has $a=b=0$ with holomorphic sectional curvature $c_0$. The
K\"ahler structure is
\begin{equation}\label{eq:ct}\begin{split}
g&=z g_0 +\frac{z^{m-1}}{\pF(z)}\, dz^2+
\frac{\pF(z)}{z^{m-1}}\,\theta^2,\qquad\qquad\;
\omega=z\omega_0+dz\wedge\theta,\\
Jdz &= \frac{\pF(z)}{z^{m-1}}\,\theta,\qquad
J\theta=-\frac{z^{m-1}}{\pF(z)}\,dz,\qquad\quad
d\theta=\omega_0,
\end{split}\end{equation}
which is the local form of a K\"ahler metric of Calabi's type~\cite{calabi0}
on a line bundle over a K\"ahler manifold $(S,g_0,\omega_0)$ of dimension
$2(m-1)$---Calabi's extremal K\"ahler metrics~\cite{calabi1} were constructed
in this way.

We remark that any such metric~\eqref{eq:ct} is conformal to a K\"ahler metric
of the same form, but with an oppositely oriented complex structure.  Indeed
we may set $\tilde z=1/z$, $\tilde g=g/z^2$ and $\tilde \pF(\tilde z)/\tilde
z^m =\pF(z)/z^m$.

\subsection{Strongly conformally Einstein K\"ahler metrics}

We return briefly to the strongly conformally Einstein manifolds of
Derdzi\'nski--Maschler~\cite{derd-masch} and explain how their classification
can be derived in a natural way in our framework.

We recall that if $(M,g,J,\omega)$ strongly conformally Einstein with
conformal factor $\tau$, then $\chi=(a\tau+b)d\tau\wedge d^c\tau/|d\tau|^2$ is
a hamiltonian $2$-form with trace $z=a\tau+b$. If $a=0$ then $\chi$ is
parallel, and $M$ is a K\"ahler product of a Riemann surface and a K\"ahler
manifold of dimension $2m-2$. We shall concentrate on the more interesting
case $a\neq0$. We set $a=1/q$ and $b=-p/q$ so that $\tau=qz+p$.  This allows
us to take the limit $q\to0$, when $\tau$ is constant and $(M,g,J,\omega)$ is
K\"ahler--Einstein. Since $\chi$ is hamiltonian of order one, with all
constant roots zero, the K\"ahler structure has the explicit
form~\eqref{eq:ct} given in the previous subsection.  This will be
conformally Einstein with conformal factor $qz+p$ if and only if
\begin{equation}\label{tosolve}
q dd^cz = \xi(z) dd^c\kappa+\eta(z) \omega,
\end{equation}
where $\kappa=\kappa_0-\frac12\log|\pF(z)|$, $\kappa_0$ is a Ricci potential
for $\omega_0$, $\xi(z)=-(qz+p)/(m-1)$ and $\eta(z)$ is an arbitrary function.
Such an equation can only hold if $\omega_0$ is K\"ahler Einstein, and we let
$dd^c \kappa_0=c\,\omega_0$ so that $\Scal_{g_0}=(m-1)c$.  We now compute
\begin{align*}
dd^cz&=
\frac{\pF}{z^{m-1}}\omega_0+\frac{z\pF'-(m-1)\pF}{z^m}
dz\wedge\theta\\ dd^c\kappa &=
\Bigl(c-\frac{\pF'}{2z^{m-1}}\Bigr)\omega_0
-\frac{z\pF''-(m-1)\pF'}{2z^m}dz\wedge\theta.
\end{align*}
Substituting in~\eqref{tosolve}, we eliminate $\eta(z)$ to obtain a
differential equation for $\pF(z)$:
\begin{multline*}
q\pF(z) - \xi(z) (cz^{m-1}-\tfrac12\pF'(z))
=\eta(z)\,z^m\\
= q(z\pF'(z)-(m-1)\pF(z)) +\tfrac12\xi(z)(z\pF''(z)-(m-1)\pF'(z)).
\end{multline*}
A particular integral is $2 c z^m/m$ (this is the flat metric), and one
solution of the homogeneous equation is $(qz - p)(qz + p)^{2m - 1}$. We can
then reduce the equation to first order to find the general homogeneous
solution as an integral:
\begin{equation*}
(qz - p)(qz + p)^{2m - 1}\Bigl(a+b\int^z
\frac{t^m}{(qt+p)^{2m}(qt-p)^2}dt\Bigr).
\end{equation*}
This is in fact a polynomial in $z$, which may be written
\begin{equation*}
\sum_{j=1}^m
\frac{j}{m}\binom{2m}{m+j}(a_+p^{m-j}q^{j-1}z^{m+j}
-a_-p^{j-1}q^{m-j}z^{m-j})
\end{equation*}
for constants $a_\pm$. The polynomial reduces to a multiple of $(qz -
p)(qz + p)^{2m - 1}$ when $q^{m+1}a_-=p^{m+1}a_+$, while for $a_+a_-=0$,
but $p,q$ nonzero, it is a hypergeometric function of $qz/p$,
essentially a Gegenbauer polynomial.

In conclusion, then, the K\"ahler structure~\eqref{eq:ct} is
conformally Einstein with conformal factor $qz+p$ if and only if
\begin{equation}
\pF(z)= \sum_{j=1}^m \frac{j}{m}\binom{2m}{m+j}
(a_+p^{m-j}q^{j-1}z^{m+j}-a_-p^{j-1}q^{m-j}z^{m-j})
+\frac{2c}m z^m
\end{equation}
for constants $a_\pm,c$ with $\rho_0=c\,\omega_0$. We have
\begin{multline*}\notag
\Scal = \frac{\Scal_{g_0}}{z}-\frac{\pF''(z)}{z^{m-1}}=
-a_+\sum_{j=1}^m \frac{j}{m}\frac{(2m)!}{(m+j-2)!(m-j)!}
p^{m-j}q^{j-1}z^{j-1}\\
+a_-\sum_{j=1}^{m-2} \frac{j}{m}\frac{(2m)!}{(m+j)!(m-j-2)!}
p^{j-1}q^{m-j}z^{-j-1},
\end{multline*}
which is polynomial in $z$ if and only if $a_-=0$, $q=0$ or $m=2$.  It is
a Killing potential (i.e., affine in $z$) if and only if it vanishes or
$q=0$ or $m=2$.

The function $Q(z)$ of Derdzi\'nski--Maschler is $\pF(z)/z^{m-1}$, and writing
$\tau=qz+p$, we recover their solution, except that we have chosen the basis
for the solutions in a way that unifies the cases $p=0$ and $p\neq 0$. We also
recall from the previous subsection that $g$ is conformal to a K\"ahler metric
$\tilde g$ with the opposite orientation.  This new K\"ahler metric is
conformally Einstein with conformal factor $\tilde\tau=q+p\tilde
z$. Just as $q=0$ corresponds to the case that $g$ is K\"ahler--Einstein, so
$p=0$ corresponds to the case that $\tilde g$ is K\"ahler--Einstein.

\appendix

\section{Conformal Killing forms and hamiltonian $2$-forms}

We have noted already that for any hamiltonian $2$-form $\phi$,
$A=\phi+\sigma\omega$ is closed. One also easily observes that
$\phi-m\sigma\omega$ is co-closed (divergence-free). In this appendix we
relate hamiltonian $2$-forms to conformal Killing $2$-forms, which have been
investigated recently by Moroianu and Semmelmann~\cite{MS,sem}.

\begin{defn} A \emph{conformal Killing} or \emph{twistor $2$-form} on a
riemannian manifold of dimension $n\geq 3$ is a $2$-form $\psi$ satisfying the
equation
\begin{equation}\label{eq:confK}
\nabla_X\psi = \frac1{n-1}X\wedge\alpha + \frac13\iota_X \beta
\end{equation}
for all vector fields $X$, where $\alpha$ is a $1$-form, $\beta$ is
a $3$-form and $\iota_X$ denotes contraction.
\end{defn}
It follows immediately from this equation that $\alpha=-\delta\psi:=\sum_i
\iota_{e_i}\nabla_{e_i}\psi$ and that $\beta=d\psi$. Hence $\psi$ is a
conformal Killing $2$-form if and only if $\nabla\phi$, which is a section of
$\Lambda^1M\otimes \Lambda^2M$, is in the kernel of the natural projections to
$\Lambda^1M$ and $\Lambda^3M$.

The notion of conformal Killing form can be extended to $p$-forms and is
conformally invariant for $p$-forms of weight $1$, i.e., for sections of
$L^{p+1}\otimes \Lambda^p M$, where $L$ is the density line bundle, i.e.,
$L^{-n}=|\Lambda^n M|$.

The following observation is due in part to Sekizawa~\cite{seki} and
Semmelmann~\cite{sem}. The essential new ingredient is~\eqref{eq:dpsi}.

\begin{prop} Let $(M,g,J,\omega)$ be a K\"ahler manifold of dimension
$n=2m\geq4$.
\begin{numlist}
\item If $\phi$ is a hamiltonian $2$-form, then the $J$-invariant
$2$-form $\psi=\phi-\frac12(\trace\phi)\omega$ is a conformal Killing
$2$-form, with $\trace\psi=\bigl(1-\frac m2\bigr)\trace\phi$.

\item Conversely, if $\psi$ is a $J$-invariant conformal Killing $2$-form,
then
\begin{equation}\label{eq:dpsi}
d\psi=-\frac3{n-1} \omega\wedge J\delta\psi
\end{equation}
and
\begin{equation}\label{eq:delpsi}
(m-2)J\delta \psi = -(2m-1)d \trace\psi.
\end{equation}
Hence
\begin{equation*}
\nabla_X\psi = -\frac1{2m-1}(X\wedge\delta\psi + JX\wedge J\delta\psi
-\g{J\delta\psi,X}\omega)
\end{equation*}
and $J\delta\psi$ is closed if $m>2$, while $\trace\psi$ is constant if
$m=2$.

\item If $\psi$ is a $J$-invariant conformal Killing $2$-form and
$J\delta\psi=df$ then $\phi=\psi+\frac1{2m-1} f\,\omega$ is a hamiltonian
$2$-form.
\end{numlist}
In particular, for $m>2$ the map $\phi\mapsto\phi-\frac12(\trace\phi)\omega$
is a bijection from hamiltonian $2$-forms to $J$-invariant conformal Killing
$2$-forms, with inverse $\psi\mapsto\psi-\frac1{m-2}(\trace\psi)\omega$.
\end{prop}
\begin{proof}
\begin{numlist}
\item If $\phi$ is hamiltonian and $\psi=\phi-\frac12\sigma\omega$, then
\begin{equation*}\notag
\nabla_X\psi=\frac12\bigl(d \sigma \wedge JX - d ^c \sigma \wedge X-
d\sigma(X)\omega\bigr)
=\frac12\bigl(X\wedge d^c\sigma-\iota_X(\omega\wedge d\sigma)\bigr).
\end{equation*}
Hence $\psi$ is a conformal Killing $2$-form with $\delta\psi=-\frac{n-1}{2}
d^c\sigma$ and $d\psi=-\frac32\omega\wedge d\sigma$.

\item Observe that $\nabla_X\psi$ is $J$-invariant and so
\begin{multline*}
\frac13 (d\psi(X,JY,Z)+d\psi(X,Y,JZ))
=\frac1{n-1}\bigl(X\wedge\delta\psi(JY,Z)+X\wedge\delta\psi(Y,JZ)\bigr)\\
=-\frac1{n-1}\bigl(\omega\wedge J\delta\psi(X,JY,Z)
+\omega\wedge J\delta\psi(X,Y,JZ)\bigr).
\end{multline*}
It follows that $\frac13 d\psi+\frac1{n-1}\omega\wedge J\delta\psi$
is a (real) $J$-invariant $3$-form, and so it must be zero---for instance
one can use the identity
\begin{align*}
2\beta(JX,Y,Z)&=\beta(Z,JX,Y)+\beta(Z,X,JY)\\
&-\beta(X,JY,Z)-\beta(X,Y,JZ)\\
&+\beta(Y,JZ,X)+\beta(Y,Z,JX).
\end{align*}
Equation~\eqref{eq:dpsi} follows immediately.
Next, the defining equation~\eqref{eq:confK} implies
\begin{equation*}\notag
(n-1)(\nabla_X\psi(JY)+\nabla_Y\psi(JX))=
\delta\psi(X)JY+\delta\psi(Y)JX-2\g{X,Y}J\delta\psi.
\end{equation*}
Taking the trace of this formula as a function of $Y$ gives
\begin{equation*}
(n-1)(2d\trace\psi(X)+J\delta\psi(X))=3J\delta\psi(X)
\end{equation*}
which is manifestly equivalent to~\eqref{eq:delpsi}.

\item This is a simple verification using~\eqref{eq:dpsi}:
\begin{align*}
\nabla_X\psi-\frac1{n-1} df(X)\,\omega
&=\frac1{n-1}\bigl(-X\wedge Jdf + \iota_X (\omega\wedge df)
-df(X)\omega\bigr)\\
&=-\frac1{n-1}\bigl(df \wedge JX - d^c f \wedge X).
\end{align*}
Note that~\eqref{eq:delpsi} gives $d\trace\phi=d\trace\psi-\frac
m{n-1}df=-\frac2{n-1}df$.
\endnumlproof

We remark that there is also a connection between hamiltonian $2$-forms and
Killing tensors (cf.~\cite{jelonek}), i.e., symmetric $2$-tensors $S$
satisfying $\sym\nabla S=0$. Indeed $\phi$ is hamiltonian if and only if
$A=\phi+\sigma\omega$ is closed and $S=J(\phi-\sigma\omega)$ is a Killing
tensor.  For the reverse implication, observe that $dA=0$ determines
$\nabla_X\phi(Y,Z)+\nabla_Y\phi(X,Z)+\nabla_Z\phi(X,Y)$ in terms of the
$1$-form $d\sigma$, while $\sym\nabla S=0$ determines
$\nabla_X\phi(JY,Z)+\nabla_Y\phi(JZ,X)+\nabla_Z\phi(JX,Y)$. Replacing $JY$ by
$Y$ in the second expression and adding to the first yields a formula for
$2\nabla_X\phi(Y,Z)+\nabla_Y\phi(Z,X)-\nabla_{JY}\phi(JZ,X)$ and the
$J$-invariant part (in $Y,Z$) is~\eqref{ham}.

\section{Vandermonde matrices}

\subsection{The inverse of a Vandermonde matrix}
A Vandermonde matrix is a $(m \times m)$-matrix of the form
\begin{equation*} V = V (\xi _1, \ldots \xi _m) =
\begin{bmatrix}
\xi _1 ^{m - 1} & \cdots & \xi _m ^{m - 1}\\
-\xi _1 ^{m - 2} & \cdots & -\xi _m ^{m - 2} \\
\vdots & \cdots & \vdots \\ (-1)^{m-1} & \cdots  & (-1)^{m-1}
\end{bmatrix},\end{equation*}
where the $\xi _j$'s are $m$ independent variables; the entries of the
Vandermonde matrix $V$ are thus defined by $V_{rj} = (-1)^{r-1}\xi_j ^{m-r}$.

We denote by $\sigma _r$ the elementary symmetric functions of the $\xi _j$'s,
so that
\begin{equation} \label{def}
\prod _{j = 1} ^m (t-\xi _j) = t^m - \sigma _1 t^{m-1} + \cdots
+ (-1) ^ m \sigma _m,
\end{equation}
for any $t$.  We also define $\sigma _0 = 1$.

Removing the variable $\xi _j$ (equivalently, differentiating with
respect to $\xi_j$) gives
\begin{equation}
\prod _{k  \neq j} (t-\xi _k) = t^{m - 1}  - \sigma _1 (\hat{\xi} _j) t
  ^{m -1} + \ldots + (-1) ^ {m-1} \sigma _{m-1} (\hat{\xi} _j),
\end{equation}
where the $\sigma _r (\smash{\hat{\xi}_j})$ are the elementary symmetric
functions of the $m-1$ variables $\xi _1, \ldots \smash{\hat{\xi}_j}, \ldots
\xi _m$ ($\xi_j$ deleted). By putting $t = \xi _i$ in the above identity, we
get
\begin{equation}
\xi _i ^{m - 1} - \sigma _1 (\hat{\xi} _j) \xi _i ^{m - 2} + \cdots
+ (- 1) ^{m - 1} \sigma _{m - 1} (\hat{\xi} _j) = \Delta _j \, \delta _{ij},
\end{equation}
where $\delta _{ij}$ is the Kronecker symbol and
$\Delta _j = \prod _{k \neq j} (\xi _j - \xi _ k)$.

This means that the matrix $W$ whose entries are
$W _{ir} = \sigma_{r-1}(\hat{\xi}_i)/\Delta_i$, i.e.,
\begin{equation*}
W  = \begin{bmatrix}
\dfrac{1}{\Delta _1} & \dfrac{\sigma _1 (\hat{\xi} _1)}{\Delta _1} & \cdots
 & \dfrac{\sigma _{m - 1} (\hat{\xi} _1)}{\Delta _1}\\
\vdots & \vdots & \cdots & \vdots \\
\dfrac{1}{\Delta _m} & \dfrac{\sigma _1 (\hat{\xi} _m)}{\Delta _m} & \cdots
 & \dfrac{\sigma _{m - 1} (\hat{\xi} _m)}{\Delta _m}
\end{bmatrix},
\end{equation*}
is a left-inverse of $V$:
\begin{equation*} \sum _{r = 1} ^m W _{ir} V _{rj} = \delta _{ij}.
\end{equation*}

\subsection{The determinant of a Vandermonde matrix}

In order to compute the determinant ${\rm det} V$ of $V$, we use the fact that
$W_{11} =1/\Delta_1$ is equal to the determinant of the minor of
$V_{11}$ in $V$ divided by $\det V$. Now the minor of $V_{11}$ is clearly $-V
(\xi _2, \ldots, \xi _m)$; we thus get the following induction formula:
\begin{equation*}
\det V (\xi _1, \ldots, \xi _m) = (-1)^{m-1}(\xi _1-\xi_2)\ldots(\xi _1-\xi _m)
\,\det V (\xi_2, \ldots, \xi _m),
\end{equation*}
from which we readily infer that
\begin{equation} \label{det}
\det V = (-1) ^{m (m - 1)/2} \prod _{i < j} (\xi _i - \xi _j).
\end{equation}

Notice that we also have
\begin{equation} \label{detdelta}
(\det V) ^2 = (-1) ^{m (m - 1)/2} \prod _{j = 1} ^m \Delta_j.
\end{equation}
Indeed, both sides are products of elements of the form $\xi_i - \xi_j$: for
each $i < j$, we get $(\xi _i - \xi _j) ^2$ in the left hand side, and $(\xi
_i - \xi _j) (\xi _j - \xi _i) = - (\xi _i - \xi _j) ^2$ in the right hand
side.

\subsection{Vandermonde identities}
In the ring of $m \times m$ matrices, a left inverse is also a
right-inverse, so that:
\begin{equation}
\sum _{j = 1} ^m V _{sj} W _{jr} = \delta _{rs};
\end{equation}
we thus obtain the following {\it Vandermonde identity}
\begin{equation} \label{id}
\sum _{j = 1} ^m \frac{(-1 ) ^{s - 1} \xi _j ^{m - s}\sigma _{r - 1}
(\hat{\xi} _j) }{\Delta _j} = \delta _{rs},
\end{equation}
for any pair $r,s=1,\ldots m$. In particular, with $r=1$, we have that
\begin{equation} \label{eq:id0}
\sum _{j = 1} ^m \frac{\xi _j ^{m - s}}{\Delta _j} = \delta_{s1}
\end{equation}
for $s=1,\ldots m$. This identity for $s=1$ may be extended to all
$s\leq 1$ to give
\begin{equation} \label{eq:id1}
\sum _{j = 1} ^m \frac{\xi _j ^{m - 1 + p}}{\Delta _j} = h_p
\end{equation}
for all $p\geq 0$, where $h_p$ is the $p$th complete symmetric function of
$\xi_1,\ldots \xi_m$. By multiplying by $t^p$, for a formal variable $t$,
and summing over $p\geq0$, this equation may be rewritten
\begin{equation*}
\sum _{j = 1} ^m \frac{\xi_j^{m-1}}{(1-\xi_j t)\Delta _j}
=\prod_{k=1}^m\frac{1}{1-\xi_k t}
\end{equation*}
where the right hand side denotes the (formal) product of
geometric series. Hence, to prove~\eqref{eq:id1}, it suffices to observe that
\begin{equation*}
\sum_{j=1}^m \xi_j^{m-1} \prod_{k\neq j} \frac{1-\xi_k t}{\xi_j-\xi_k} = 1.
\end{equation*}
This follows because the left hand side is a polynomial in $t$, of degree at
most $m-1$, whose value at $t=1/\xi_j$ is equal to $1$ for all $j=1,\ldots m$.
(In fact, this is more or less the Lagrange interpolation formula.)

Similarly, we can extend~\eqref{id} to obtain
\begin{equation}\label{eq:id2}
\sum _{j = 1} ^m \frac{\xi _j ^{m + k}\sigma_{r-1}(\hat\xi_j)}{\Delta _j}
= \sum_{s=0}^k (-1)^s h_{k-s}\sigma_{r+s}
\end{equation}
for all $r=1,\ldots m$ and all $k\geq 0$. Here, by convention,
$\sigma_{r+s}=0$ for $r+s>m$. We reduce~\eqref{eq:id2}
to~\eqref{eq:id1} by means of the obvious identity:
\begin{equation*}
\xi_j^{m+k} \sigma_{r-1}(\hat \xi_j)
= \sum_{s=0}^{m-r} (-1)^s \xi_j^{m-1+k-s} \sigma_{r+s}.
\end{equation*}
(Evidently $\xi_j \sigma_{r-1}(\hat \xi_j) = \sigma_{r}- \sigma_{r}(\hat
\xi_j)$.)  Substitute this into the left hand side of~\eqref{eq:id2}, and note
that the summation over $s$ can be made from $0$ to $k$, using the Vandermonde
identity~\eqref{eq:id0} to eliminate any extra terms. Now
applying~\eqref{eq:id1} for each $s$ yields the right hand side
of~\eqref{eq:id2}.

There is one further identity we shall need, namely
\begin{equation}\label{eq:dxi}
\frac{\partial}{\partial\xi_i}\Bigl(
\sum _{j = 1}^m \frac{\xi _j ^{m + k}\sigma_{r-1}(\hat\xi_j)}{\Delta _j}\Bigr)
=\sigma_{r-1}(\hat\xi_i)\sum_{s=0}^k h_{k-s}\xi_i^s.
\end{equation}
We prove this using~\eqref{eq:id2}: multiplying by $t^k$ and summing
over $k$, it suffices to show
\begin{equation*}
\sum_{s\geq 0} \frac{\partial}{\partial\xi_i}
\Bigl(\frac{  (-1)^s \sigma_{r+s} t^s }
{ \prod_{j=1}^m (1-\xi_j t) }\Bigr)
= \frac{  \sigma_{r-1}(\hat \xi_i)}
{ (1-\xi_i t) \prod_{j=1}^m (1-\xi_j t)}.
\end{equation*}
This holds since direct computation of the left hand side gives
\begin{equation*}
\sum_{s\geq 0} \frac{(-1)^s \sigma_{r+s-1}(\hat \xi_i) t^s
+(-1)^s\sigma_{r+s}(\hat\xi_i) t^{s+1}}
{ (1-\xi_i t) \prod_{j=1}^m (1-\xi_j t)}
\end{equation*}
using $\sigma_{r+s}(\hat\xi_i)=\sigma_{r+s}-\xi_i\sigma_{r+s-1}(\hat\xi_i)$.
All terms now cancel in pairs except the first one with $s=0$.

In fact we shall only make serious use of the identities~\eqref{eq:id2}
and~\eqref{eq:dxi} for $0\leq k\leq 2$. In particular~\eqref{eq:id2} implies
\begin{align}
\sum _{j = 1} ^m \frac{\xi ^m _j}{\Delta _j} \,\sigma _{r - 1} (\hat{\xi} _j)
&= \sigma _r,\\
\sum_{j=1}^m \frac{\xi_j^{m+1}}{\Delta _j} \, \sigma _{r - 1} (\hat{\xi} _j)
&= \sigma _1 \sigma _r - \sigma _{r + 1},\\
\sum_{j=1}^m \frac{\xi_j^{m+2}}{\Delta _j} \, \sigma _{r - 1} (\hat{\xi} _j)
&= (\sigma_1^2-\sigma_2)\sigma _r-\sigma_1 \sigma_{r + 1}+\sigma _{r+2}.
\end{align}


\begin{thebibliography}{99}

\bibitem{Abreu} M. Abreu,
{\it K{\"a}hler geometry of toric varieties and extremal metrics},
Int. J. Math. {\bf 9}  (1998) 641--651.

\bibitem{AG2} V. Apostolov and P. Gauduchon, {\it Self-dual Einstein Hermitian
$4$-manifolds}, Ann. Scuola Norm. Sup. Pisa (5) {\bf 1} (2002) 203--243.

\bibitem{ACG0} V. Apostolov, D.~M.~J. Calderbank and P. Gauduchon,
{\it The geometry of weakly self-dual K\"ahler surfaces}, Compositio
Math. {\bf 135} (2003), 279--322.

\bibitem{besse} A. L. Besse, {\it Einstein manifolds},
Ergeb. Math. Grenzgeb. {\bf 3}, Springer-Verlag, Berlin, Heidelberg,
New York, 1987.

\bibitem{bryant} R. Bryant, {\it Bochner--K{\"a}hler metrics},
J. Amer. Math. Soc. {\bf 14} (2001) 623--715.

\bibitem{calabi0} E. Calabi, {\it M\'etriques k\"ahl\'eriennes et
fibr\'es holomorphes}, Ann. Sci. Ecole Norm. Sup. (4) {\bf 12} (1979) 269--294.

\bibitem{calabi1} E. Calabi, {\it Extremal K{\"a}hler metrics},
Seminar on Differential Geometry, Princeton University Press 1982.

\bibitem{derd-masch} A. Derdzi\'nski and G.~Maschler, {\it Local types of
conformally-Einstein K\"ahler metrics in higher dimensions}, Preprint (2000),
available at arXiv:math.DG/0204013.

\bibitem{DH} J. J. Duistermaat and G. J.  Heckman, {\it On the variation in
the cohomology of the symplectic form of the reduced phase space},
Invent. Math. {\bf 69} (1982) 259--268.

\bibitem{Gray} A.~Gray,
{\it Pseudo-Riemannian almost product manifolds and submersions},
J. Math. Mech. {\bf 16} (1967) 715--737.

\bibitem{G:kstv} V.~Guillemin,
{\it K{\"a}hler structures on toric varieties},
J. Diff. Geom.  {\bf 40}  (1994)  285--309.

\bibitem{hwang-singer} A.~D. Hwang and M.~A. Singer,
{\it A momentum construction for circle-invariant Kahler metrics},
Trans. Amer. Math. Soc. {\bf 354} (2002) 2285--2325.

\bibitem{jelonek} W. Jelonek,
{\it Compact K{\"a}hler surfaces with harmonic anti-self-dual Weyl tensor},
Diff. Geom. Appl. 16 (2002), 267--276.

\bibitem{koisak} N. Koiso and Y. Sakane, {\it Nonhomogeneous
K\"ahler--Einstein metrics on compact complex manifolds}, in
Curvature and Topology of Riemannian Manifolds (Kataka, 1985),
Lecture Notes in Math. 1201, Springer-Verlag, Berlin, 1986, 165--179.

\bibitem{kostant} B. Kostant, {\it Holonomy and the Lie algebra of
infinitesimal motions of a Riemann manifold},
Trans. Amer. Math. Soc. {\bf 80} (1955) 528--542.

\bibitem{lebrun} C.~R. LeBrun, {\it Explicit self-dual metrics on
{$\C P^2\#\cdots\#\C P^2$}}, J.~Diff. Geom. {\bf 34} (1991) 223--253.

\bibitem{mats-tanno} M. Matsumoto and S. Tanno,
{\it On K{\"a}hler spaces with parallel or vanishing Bochner curvature
tensor}, Tensor N. S. {\bf 27} (1973) 291--294.

\bibitem{MS} A. Moroianu and U. Semmelmann,
{\it Twistor forms on K\"ahler manifolds},
available at arXiv:math.DG/0204322

\bibitem{ONeill} B.~O'Neill,
{\it The fundamental equations of a submersion},
Michigan Math. J. {\bf 13} (1966) 459--469.

\bibitem{ped-poon} H. Pedersen and Y.~S. Poon,
{\it Hamiltonian construction of K\"ahler--Einstein metrics and
K\"ahler metrics of constant scalar curvature},
Comm. Math. Phys. {\bf 136} (1991) 309--326.

\bibitem{seki} M.~Sekizawa,
{\it On conformal Killing tensors of degree $2$ in K\"ahlerian spaces},
TRU Math. {\bf 6} (1970) 1--5.

\bibitem{sem} U.~Semmelmann, {\it Conformal Killing forms on Riemannian
manifolds}, Habilitationschrift, Universit\"at M\"unchen (2002), available at
arXiv:math.DG/0206117

\bibitem{christina} C. T{\o}nnesen-Friedman,
{\it Extremal K{\"a}hler metrics on minimal ruled surfaces}, J. reine
angew. Math. {\bf 502} (1998) 175--197.

\bibitem{christina2} C. T{\o}nnesen-Friedman,
{\it Extremal Kahler metrics and Hamiltonian functions} II,
Glasg. Math. J. {\bf 44} (2002) 241--253.

\bibitem{wang} M.~Y. Wang,
{\it Einstein metrics from symmetry and bundle constructions},
Surv. Diff. Geom., VI,
Int. Press, Boston, MA, 1999, 287--325.

\end{thebibliography}
\end{document}